\documentclass[a4]{amsart}
\usepackage{amsthm,amssymb,amsmath,color}
\usepackage{pdfsync}
\newcommand{\blue}{\color{darkblue}}

\newcommand{\levy}{L\'evy}

\newcommand{\clt}{central limit theorem}

\newcommand{\garch}{{\rm GARCH}$(1,1)$}

\newcommand{\sta}{St\u aric\u a}
\newcommand{\ex}{{\rm e}\,}

\newcommand{\asy}{asymptotic}

\newcommand{\ts}{time series}

\newcommand{\ul}{\underline}

\definecolor{darkblue}{rgb}{.1, 0.1,.8}
\definecolor{darkgreen}{rgb}{0,0.8,0.2}
\definecolor{darkred}{rgb}{.8, .1,.1}

\newtheorem{lemma}{Lemma}[section]

\newtheorem{theorem}[lemma]{Theorem}

\newcommand{\bbc}{{\mathbb C}}
\newtheorem{proposition}[lemma]{Proposition}
\newtheorem{definition}[lemma]{Definition}
\newtheorem{corollary}[lemma]{Corollary}
\newtheorem{example}[lemma]{Example}
\newtheorem{exercise}[lemma]{Exercise}
\newtheorem{remark}[lemma]{Remark}
\newtheorem{fig}[lemma]{Figure}
\newtheorem{tab}[lemma]{Table}

\newcommand{\MC}{Markov chain}

\newcommand{\bth}{\begin{theorem}}
\newcommand{\ethe}{\end{theorem}}

\newcommand{\bre}{\begin{remark}\em }
\newcommand{\ere}{\end{remark}}

\newcommand{\ble}{\begin{lemma}}
\newcommand{\ele}{\end{lemma}}
\newcommand{\sre}{stochastic recurrence equation}
\newcommand{\pp}{point process}
\newcommand{\bde}{\begin{definition}}
\newcommand{\ede}{\end{definition}}
\newcommand{\bco}{\begin{corollary}}
\newcommand{\eco}{\end{corollary}}

\newcommand{\bpr}{\begin{proposition}}
\newcommand{\epr}{\end{proposition}}

\newcommand{\bexer}{\begin{exercise}}
\newcommand{\eexer}{\end{exercise}}

\newcommand{\bexam}{\begin{example}}
\newcommand{\eexam}{\end{example}}

\newcommand{\bfi}{\begin{fig}}
\newcommand{\efi}{\end{fig}}

\newcommand{\btab}{\begin{tab}}
\newcommand{\etab}{\end{tab}}

\newcommand{\fidi}{finite-dimensional distribution}
\newcommand{\rv}{random variable}

\newcommand{\rhs}{right-hand side}
\newcommand{\df}{distribution function}

\newcommand{\beao}{\begin{eqnarray*}}
\newcommand{\eeao}{\end{eqnarray*}\noindent}

\newcommand{\beam}{\begin{eqnarray}}
\newcommand{\eeam}{\end{eqnarray}\noindent}

\newcommand{\beqq}{\begin{equation}}
\newcommand{\eeqq}{\end{equation}\noindent}

\newcommand{\bce}{\begin{center}}
\newcommand{\ece}{\end{center}}

\newcommand{\barr}{\begin{array}}
\newcommand{\earr}{\end{array}}

\newcommand{\stp}{\stackrel{P}{\rightarrow}}
\newcommand{\std}{\stackrel{d}{\rightarrow}}

\newcommand{\stv}{\stackrel{v}{\rightarrow}}
\newcommand{\stw}{\stackrel{w}{\rightarrow}}

\newcommand{\eqd}{\stackrel{d}{=}}

\newcommand{\vague}{\stackrel{\lower0.2ex\hbox{$\scriptscriptstyle
                    \it{v} $}}{\rightarrow}}
\newcommand{\weak}{\stackrel{\lower0.2ex\hbox{$\scriptscriptstyle
                    \it{w} $}}{\rightarrow}}
\newcommand{\what}{\stackrel{\lower0.2ex\hbox{$\scriptscriptstyle
                    \it{\hat{w}} $}}{\rightarrow}}

\newcommand{\bdis}{\begin{displaymath}}
\newcommand{\edis}{\end{displaymath}\noindent}

\newcommand{\N}{\mathbb{N}}
\newcommand{\R}{\mathbb{R}}

\newcommand{\nto}{n\to\infty}
\newcommand{\kto}{k\to\infty}
\newcommand{\xto}{x\to\infty}

\newcommand{\ov}{\overline}
\newcommand{\un}{\underline}
\newcommand{\wt}{\widetilde}

\newcommand{\vep}{\varepsilon}

\newcommand{\la}{\lambda}

\newcommand{\regvary}{regularly varying}
\newcommand{\slvary}{slowly varying}
\newcommand{\regvar}{regular variation}

\newcommand{\bbr}{{\mathbb R}}

\newcommand{\bbz}{{\mathbb Z}}
\newcommand{\bbn}{{\mathbb N}}

\newcommand{\bbs}{{\mathbb S}}

\newcommand{\con}{convergence}

\newcommand{\evt}{extreme value theory}

\newcommand{\st}{such that}
\newcommand{\fif}{if and only if}
\newcommand{\wrt}{with respect to}
\newcommand{\chf}{characteristic function}
\newcommand{\fct}{function}

\newcommand{\ds}{distribution}

\newcommand{\rep}{representation}

\newcommand{\seq}{sequence}

\newcommand{\ins}{insurance}

\newcommand{\pro}{probabilit}

\newcommand{\ms}{measure}

\newcommand{\ld}{large deviation}
\newcommand{\bfx}{{\bf x}}

\newcommand{\bfy}{{\bf y}}

\newcommand{\bfs}{{\bf s}}

\def\1{\ensuremath{\mathrm{1}\hspace{-.35em} \mathrm{1}}} 

\def\E{{\mathbb E}}

\def\N{\mathbb{N}}
\def\P{{\mathbb{P}}}
\def\R{\mathbb{R}}
\def\Z{\mathbb{Z}}

\renewcommand{\le}{\ensuremath{\leqslant}}
\renewcommand{\ge}{\ensuremath{\geqslant}}

\newcommand{\introo}[2]{{\left]{#1,\,#2\,}\right[\kern1pt}}

\newcommand{\intrfo}[2]{{\left[{#1,\,#2}\right[\kern1pt}}

\begin{document}

\title[A large deviations approach to limit theory for heavy-tailed
time series]{A large deviations approach to limit theory for heavy-tailed
time series}
\thanks{
Both authors would like to thank their home institutions for
hospitality when visiting each other. Parts of this paper were written
when T. Mikosch was on sabbatical in 2013. He would like to thank
the Columbia Statistics Department and the Forschungsinstitut f\"ur
Mathematik at ETH Z\"urich for their kind hospitality. The paper was
finished when T. Mikosch visited Universit\'e Marie and Pierre Curie
and Institut H. Poincar\'e. He would like to thank Paul Doukhan for
providing generous financial support. Both authors would like to thank Paul
Doukhan and Adam Jakubowski for numerous discussions on the topic of
this paper. Financial supports by the Danish Organization for Free Research
(DFF) Grant 107959 and by the ANR network AMERISKA are gratefully acknowledged.}

\author[Thomas Mikosch and Olivier Wintenberger]{T. Mikosch and O. Wintenberger}
\address{T. Mikosch\\ 
University of Copenhagen\\Universitetsparken 5\\
2100 Copenhagen\\Denmark}
\email{mikosch@math.ku.dk}
\address{O. Wintenberger\\Sorbonne Universit\'es\\ 
UPMC Univ Paris 06\\ LSTA, Case 158
4 place Jussieu\\
75005 Paris\\ France}
\email{olivier.wintenberger@upmc.fr}

\begin{abstract}
In this paper we propagate a \ld s approach for proving limit
theory for (generally) multivariate time series with heavy tails. We make
this notion precise by introducing \regvary\ time series.
We provide general \ld\ results for \fct als acting on a sample path
and
vanishing in some neighborhood of the origin. 
We study a variety
of such \fct als, including large deviations of 
random walks, their suprema, the ruin \fct al, and further 
derive weak limit theory for maxima, \pp es, cluster \fct als and the tail
empirical process. One of the main results of this paper concerns
bounds for the ruin probability in various heavy-tailed models including GARCH, stochastic volatility models and solutions to \sre s.
\end{abstract}


 \maketitle

\section{Preliminaries and basic
motivation}\label{sec:intr}\setcounter{equation}{0}
In the last decades, a lot of efforts has been put into the
understanding of limit theory for dependent \seq s, including \MC s
(Meyn and Tweedie \cite{meyn:tweedie:1993}), weakly dependent
\seq s (Dedecker et al. \cite{dedecker:doukhan:lang:leon:louhichi:prieur:2007}),
long-range dependent
\seq s (Doukhan et al. \cite{doukhan:oppenheim:taqqu:2003},
Samorodnitsky \cite{samorodnitsky:2006}),
empirical processes (Dehling et al. \cite{dehling:mikosch:sorensen:2002})
and more general structures (Eberlein and Taqqu \cite{eberlein:taqqu:1986}), to name a few
references. A smaller part of the theory was devoted to limit theory
under extremal dependence for \pp es, maxima, partial sums, tail
empirical processes. Resnick \cite{resnick:1986,resnick:1987}
started a systematic study of the relations between the \con\ of
\pp es, sums and maxima; see also Resnick \cite{resnick:2007}
for a recent account. He advocated the use of {\em multivariate
  \regvar } as a flexible tool to describe heavy-tail phenomena
combined with advanced continuous mapping techniques. For example,
maxima and sums are understood as \fct als acting on an underlying
\pp ; if the \pp\ converges these \fct als converge as well and their limits
are described in terms of the points of the limiting \pp .
\par
Davis and Hsing \cite{davis:hsing:1995} recognized the power
of this approach for limit theory of \pp es, maxima, sums, and large
deviations for dependent {\em \regvary } processes, i.e., stationary
\seq s whose \fidi s are \regvary\ with the same index.
Before \cite{davis:hsing:1995}, limit theory for particular
\regvary\ stationary \seq s was
studied for the sample mean,
maxima, sample autocovariance and autocorrelation
\fct s of linear and bilinear processes with iid \regvary\ noise 
and \evt\ was considered for \regvary\ ARCH processes and solutions to
\sre ;
see Rootz\'en
\cite{rootzen:1976}, 
Davis and Resnick
\cite{davis:resnick:1985,davis:resnick:1985a,davis:resnick:1986,davis:resnick:1996}, de
Haan et al. \cite{haan:resnick:rootzen:vries:1989}. Davis and Hsing
\cite{davis:hsing:1995} introduced \regvar\ of a random \seq\
as a general flexible tool for proving limit theory
for heavy-tail phenomena. The theory in
\cite{davis:hsing:1995} is a benchmark for the results of the present
paper. We quote the main result of \cite{davis:hsing:1995} on convergence of point processes
for reasons of comparison (see Theorem~\ref{thm:davhs} below); 
we will also give a 
short alternative proof in this paper. The main result of
  \cite{davis:hsing:1995} has been used to derive the \clt\ with
  infinite variance stable limit  via the mapping theorem. When
  studying other functionals of the sample path than sums and maximas,
  this approach is limited by the continuity condition in the mapping
  theorem. For example, 
the asymptotics for the supremum of the partial sums  
follows only under additional restrictions from the point process
approach;  see Basrak et al. \cite{basrak:krizmanic:segers:2012}.
\par
We introduce a new approach to bypass this restriction. We essentially follow an argument of Ja\-ku\-bows\-ki 
\cite{jakubowski:1993,jakubowski:1997} and 
 Jakubowski and Kobus \cite{jakubowski:kobus:1989},
using a telescoping sum
approach. Under suitable anti-clustering conditions, 
this argument can be applied to Laplace \fct als, \chf s of
sums,  \df s of maxima, etc. 
This approach turned out to be fruitful in our previous work;
see Bartkiewicz et al.
\cite{bartkiewicz:jakubowski:mikosch:wintenberger:2011},
Mikosch and Wintenberger
\cite{mikosch:wintenberger:2013,mikosch:wintenberger:2012}.  A careful
study of related work such as Davis and Hsing \cite{davis:hsing:1995},
Basrak and Segers \cite{basrak:segers:2009}, 
Segers \cite{segers:2005}, 
Balan and Louhichi \cite{balan:louhichi:2009,balan:louhichi:2010} and  Yun \cite{yun:2000},
shows that the telescoping approach has been used in these papers as well.
The aim of this paper is to understand the common structural 
properties of these results and their close relationship with large
deviation theory.
\par
The framework of this paper is the one of regularly varying stationary processes that we introduce now. We commence with a random vector $X$ with values in $\R^d$ for some $d\ge 1$.
We say that this vector (and its \ds ) are {\em \regvary\ with index
$\alpha >0$} if the following relation holds as $\xto$:
\begin{equation}\label{eq:rv}
\frac{\P(|X| > ux, X/| X| \in\cdot )}{\P(|X| > x)}\stw
u^{-\alpha}\, \P(\Theta\in\cdot),\quad u>0\,.
\end{equation}
Here $\stw$ denotes weak \con\ of finite \ms s  and $\Theta$ is a vector with values in the unit sphere
$\bbs^{d-1} = \{x\in\bbr^d : |x| = 1\}$ of   $\R^d$. Its \ds\ is the {\em spectral \ms }
of \regvar\ and depends on the choice of the norm. However, the definition of \regvar\ does not
depend on any concrete norm; for convenience we always refer to the Euclidean norm. An equivalent way to define
\regvar\ of $X$ is to require that there exists a non-null Radon \ms\
$\mu_X$ on the
Borel $\sigma$-field of $\ov \bbr_0^d=
\ov \bbr^d\setminus\{\bf0\}$ \st\
\beam\label{eq:1}
n\, \P(a_n^{-1} X\in \cdot )\stv \mu_X\,,
\eeam
where the \seq\ $(a_n)$ can be chosen \st\ $n\,\P(|X|>a_n)\sim 1$ and $\stv$ refers to
vague \con . The limit \ms\ $\mu_X$ necessarily has the property
$\mu_X(u\cdot)= u^{-\alpha}\mu_X(\cdot)\,,u>0$,
which explains the relation with the index $\alpha$.
We refer to Bingham et al.
\cite{bingham:goldie:teugels:1987} for an encyclopedic treatment of
one-dimensional \regvar\ and Resnick
\cite{resnick:1987,resnick:2007} for the multivariate case.
\par
Next consider a strictly stationary \seq\
$(X_t)_{t\in\Z}$ of $\bbr^d$-valued random vectors with a generic
element $X$.
It is {\em \regvary\ with index $\alpha>0$} if every lagged vector
$(X_1, . . . ,X_k)$, $k\ge 1$, is \regvary\ in the sense of \eqref{eq:rv}; see
Davis and Hsing \cite{davis:hsing:1995}. An equivalent description of a \regvary\ \seq\ $(X_t)$ is achieved by
exploiting \eqref{eq:1}: for every $k\ge 1$, there exists a non-null Radon \ms\ $\mu_k$ on
 the Borel $\sigma$-field of $ \ov \bbr_0^{dk}$ \st\
\beam\label{eq:1a}
n\, \P(a_n^{-1} (X_1,\ldots,X_k)\in \cdot )\stv \mu_k\,,
\eeam
where $(a_n)$ is chosen \st\ $n\,\P(|X_0|>a_n)\sim 1$.
\par

A convenient characterization of a \regvary\ \seq\  $(X_t)$ was given in
Theorem 2.1 of Basrak and Segers \cite{basrak:segers:2009}: there exists a \seq\ of $\bbr^d$-valued random vectors
$(Y_t)_{t\in\bbz}$ \st\ $\P(|Y_0| > y) = y^{-\alpha}$ for $y > 1$ and
for $k\ge 0$,
\beao
\P (x^{-1}(X_{-k},\ldots,X_k)\in\cdot \mid |X_0| > x)\stw
\P((Y_{-k},\ldots,Y_k)\in \cdot)\,,\quad \xto\,.
\eeao
The process
$(Y_t)$ is the {\em tail process} of $(X_t)$. Writing $\Theta_t = Y_t/|Y_0|$ for $t\in \Z$,
one also has for $k\ge 0$,
\beam\label{eq:sepctraltailprocess}
\P( |X_0|^{-1}(X_{-k},\ldots,X_k)\in\cdot \mid |X_0| > x)
\stw \P ((\Theta_{-k},\ldots,\Theta_k)\in \cdot)\,,\quad \xto\,.
\eeam
 We  identify
$|Y_0|\,(Y_t/|Y_0|)_{|t|\le k}= |Y_0|\,(\Theta_t)_{|t|\le k}$, $k\ge
0$. Then $|Y_0|$ is independent of $(\Theta_t)_{|t|\le k}$ for every
$k\ge 0$. We refer
to $(\Theta_t)_{t\in\bbz}$ as
the {\em spectral tail process} of $(X_t)$.  In what follows, we will
make heavy use of the tail and spectral tail processes: most of our
results will be expressed in terms of them. We will refer to either
condition \eqref{eq:1a} and the equivalent tail and spectral tail
conditions  as  ${\bf (RV_\alpha)}$.
\par
The condition ${\bf (RV_\alpha)}$ is equivalent to the fact 
that for any $\vep>0$, any continuous bounded function $f(x_0,x_1,\ldots)$ 
on $(\R^d)^\N$ which vanishes for $|x_0|\le \vep$ the following
relation holds:
\beao
\frac{\E[f(x^{-1}(X_0,\ldots,X_k,0,0,\ldots))]}{\P(|X_0|>x)}\to \int_0^\infty \E[f(y\Theta_0,\ldots,y\Theta_k,0,0,\ldots)]\,d(-y^{-\alpha}),\qquad k\ge 0.
\eeao
When considering the extremal properties of a sample path $(X_0,\ldots,X_n)$ for large $n$, it
is not natural to assume that $X_0$ is large. To overcome this restriction
a telescoping argument helps. 
It shows that for any $\vep>0$, any continuous bounded function 
$f$ on $(\R^d)^{\N}$ which vanishes if $|x_t|\le \vep$ for all $t\ge
0$, 
the following relation holds:
\begin{align}
\frac{\E[f(x^{-1}(X_0,\ldots,X_n,0,0,\ldots))]}{\P(|X_0|>x)}\to\nonumber& \sum_{j=0}^n\int_0^\infty \E[f(\underbrace{0,\cdots,0}_{j},y\Theta_0,\ldots,y\Theta_{n-j},0,0,\ldots)\\&-f(\underbrace{0,\cdots,0}_{j+1},y\Theta_1,\ldots,y\Theta_{n-j},0,0,\ldots)]\,d(-y^{-\alpha}),\qquad n\ge 0\label{eq:dec3}.
\end{align}
The \rhs\ is no longer well-behaved when $\nto$. 
To study the asymptotic extremal properties of the sample path a
C\`esaro argument is required; under suitable assumptions the \rhs\ 
will converge after renormalization with $n+1$. In addition, the mean ergodic theorem can be used if we assume conditions such as
$$
f(\underbrace{0,\cdots,0}_{j},y\Theta_0,\ldots,y\Theta_{n-j},0,0,\ldots)=f(y\Theta_0,\ldots,y\Theta_{n-j},0,0,\ldots).
$$
Such a relation is satisfied under our condition {\bf (C$_\vep$)}; see Section \ref{sec:ld}. If the spectral tail
process 
is also ergodic we obtain
\beam\label{eq:dec3a}
\lim_{n\to\infty}\lim_{x\to \infty}\frac{\E[f(x^{-1}(X_0,\ldots,X_n,0,0,\ldots))]}{(n+1)\P(|X_0|>x)}= \int_0^\infty \E[f(y(\Theta_t)_{t\ge 0}) -f((y\Theta_t)_{t\ge 1})]\,d(-y^{-\alpha}).\nonumber\\
\eeam
Our large deviation approach can be seen as an equicontinuity argument
applied to the intermediate result \eqref{eq:dec3}. Under suitable
assumptions
it is possible that \eqref{eq:dec3a} 
 holds uniformly on a region $\Lambda_n$ of $x$-values when
 $n\to\infty$. 
We will use an anti-clustering condition to enforce the equicontinuity and  ergodicity of the spectral tail process.
Thanks to this new approach we can characterize 
the large deviations for various functionals of the sample path. 
For example, we obtain the large deviations of the supremum of the 
partial sums and we derive the asymptotic behavior of the ruin
probability.
\par 
This new approach describes the limiting
behavior of extremes of \regvary\ \seq s in term of their spectral tail processes. We refer to
calculations of the spectral tail processes for concrete examples such that certain Markov, stochastic volatility, \garch\ processes,
solutions to \sre s, max-stable processes, and other examples in Basrak
et al.  \cite{basrak:segers:2009,basrak:krizmanic:segers:2012},
Mikosch and Wintenberger \cite{mikosch:wintenberger:2012}, Davis et
al. \cite{davis:mikosch:zhao:2013}; see also
Examples~\ref{exam:goldie}--\ref{exam:sv} below.
\par
The paper is organized as follows.
In Section~\ref{sec:prel} we provide some probabilistic tools 
used in  the article and formulate the main result of Davis and Hsing
\cite{davis:hsing:1995}. In Section~\ref{sec:ld} we formulate our main result 
about the \ld s for \fct als acting on the sample paths of a \regvary\
\seq ; see Theorem~\ref{th:antic}.
In Section~\ref{sebsec:supremaandruin} we show two other main results
of this paper. In Theorem~\ref{thm:suprema} we give a uniform \ld\
bound for the suprema of a random walk constructed from a \regvary\
\seq . A modification of the proof of  Theorem~\ref{thm:suprema}
 is then used to give bounds for the tails of the
ruin \fct al; see Theorem~\ref{thm:ruina}. We apply the latter result
to solutions to \sre s, \garch\ and stochastic volatility  processes.
In Section~\ref{subsec:clfctal}
we show how the \ld\ approach helps to prove results for
cluster \fct als and in Section~\ref{sec:tailemp} we apply \ld s to
get results for the tail empirical process of a \regvary\
\seq .
Finally,
the proofs of the results are provided in Section~\ref{sec:proofs}.

\section{Preliminaries}\label{sec:prel}\label{sec:kdep}

\subsection{Anti-clustering conditions}
Davis and Hsing  \cite{davis:hsing:1995} introduced a condition that avoids
``long-range dependence'' of high level exceedances of the process
$(X_t)$:\\
{\em Anti-clustering condition} {\bf (AC)}:
Let $m=m_n\to\infty$ be an integer \seq\ \st\ $m_n=o(n)$ as $\nto$
and $(a_n)$ the normalizing \seq\ from ${\bf (RV_\alpha)}$. They assume
that
\beam\label{eq:accond}
\lim_{k\to\infty}\limsup _{\nto} \P\big(\wt M_{k, m_n} >\delta
a_n\mid |X_0|>\delta a_n\big)=0\,,\quad \delta>0\,, 
\eeam
where 
\beam\label{eq:maxst}
M_{s,t}=\max_{s\le i \le t} |X_i|\,,\quad s\le t\,,\qquad \mbox{and}\qquad \wt M_{s,t}=\max_{s\le|i| \le t} |X_i|\,,\quad s\le t.
\eeam
This condition assures that extremal clusters of $(X_t)$ get separated
from each other when time goes by, i.e., the influence of an extremal
shock at some time  does not last forever.
Conditions of this type are common in the extreme value literature,
e.g. the popular condition $D'(a_n)$; see
Leadbetter et al. \cite{leadbetter:lindgren:rootzen:1983}, cf. Section
4.4 in Embrechts et al. \cite{embrechts:kluppelberg:mikosch:1997}.
It is often easy to verify  {\bf (AC)} by checking
\beao
\lim_{k\to\infty}\limsup _{\nto} \sum_{k\le |t| \le m_n} \P(|X_t|>
\delta a_n\mid |X_0|>  \delta a_n)=0\,,\quad \delta>0\,.
\eeao
The following result
can be found in Segers \cite{segers:2005} and Basrak and Segers
\cite{basrak:segers:2009}; see also O'Brien \cite{obrien:1987}.
\bpr\label{prop:1} Let $(X_t)$ be a non-negative strictly stationary \seq\ which is
\regvary\ with index $\alpha>0$ and satisfies
{\bf (AC)}. Then
the limit
\beao
 \lim_{\nto} \P(M_{1,m_n}\le \delta a_n\mid X_0>\delta
 a_n)=\gamma
= P\big(\sup_{t\ge 1}Y_{t}\le 1\big)=P\big(\sup_{t\ge 1}Y_{-t}\le 1\big)\,,
\eeao
exists for every $\delta>0$, it is positive and
$\gamma$ is the extremal index of $(X_t)$.
\epr
The extremal index of a real-valued stationary \seq\ is often
interpreted as reciprocal of the expected cluster size of 
high level exceedances. This intuition can be made precise; see for
example the monographs Leadbetter et al. \cite{leadbetter:lindgren:rootzen:1983} and Embrechts et al. 
\cite{embrechts:kluppelberg:mikosch:1997}, Section 8.1.

\subsection{Mixing conditions}
Davis and Hsing \cite{davis:hsing:1995} assumed a mixing condition in
terms of Laplace \fct als of the \pp es
\beam\label{eq:nn}
N_{nj}=\sum_{t=1}^j \vep_{a_n^{-1} X_t}\,, \quad j=1,\ldots,n\,,\quad
N_{nn}=N_n\,,\quad n\ge 1\,,
\eeam
with state space $\ov \bbr_0^d=\ov
\bbr^d\setminus \{0\}$, where $\ov \bbr=\bbr\cup\{-\infty,\infty\}$.
This condition reads as follows:\\[1mm]
{\em Condition} ${\mathcal A}(a_n)$: For the same  \seq\ $(m_n)$ as in
{\bf (AC)} and with $k_n=[n/m]\to\infty$,
\beao
\E  \ex^{- \int f dN_n}- \Big(\E  \ex^{- \int f dN_{n,m_n}}\Big)^{k_n}\to
0\,,\quad f\in \bbc_K^+\,,
\eeao
where $\bbc_K^+$ is the set of non-negative continuous \fct s with
compact support.

Boundedness of a subset of $\ov\bbr_0^d$
means that it is bounded away from zero, in particular, compact sets in $\ov\bbr_0^d$
are bounded away from zero.
Davis and Hsing
  \cite{davis:hsing:1995} assumed a slightly more general version of
${\mathcal A}(a_n)$: their \fct s $f$ are any non-negative step \fct s on
$\ov\bbr_0^d$ with bounded support. For both classes of \fct s, the \con\ of the
Laplace \fct als $\E \ex^{-\int f dQ_n}\to \E \ex^{-\int f dQ}$
for \pp es $(Q_n),Q$ is equivalent to $Q_n\std Q$; see Kallenberg
\cite{kallenberg:1983}. The restriction to  $f\in \bbc_K^+$ is common
in the literature, e.g. Resnick
\cite{resnick:1986,resnick:1987,resnick:2007}. Various papers which build on
Davis and Hsing  \cite{davis:hsing:1995} also assume  ${\mathcal A}(a_n)$
for $f\in\bbc_K^+$; see Basrak and Segers \cite{basrak:segers:2009},
Balan and Louhichi \cite{balan:louhichi:2009}.
Condition ${\mathcal A}(a_n)$ for suitable \seq s $(m_n)$ follows from both 
strong mixing and  weak 
dependence in the sense of Dedecker and
Doukhan~\cite{dedecker:doukhan:2003}.
\bre\label{rem:id}
Let $B_\delta^c=\{x\in \ov\bbr^d_0: |x|>\delta\}$, $\delta>0$.
Under \regvar\ of $(X_t)$, $P(N_{nm}(B_\delta ^c)>0)\le m_n \P(| X_1|>\delta
a_n)\to 0$ and therefore an iid \seq\
$(\wt N_{nm}^{(i)})$ of copies of $N_{nm}$ is a {\em null-array} in
the sense of Kallenberg \cite{kallenberg:1983}. Then, according to
Theorem 6.1 in \cite{kallenberg:1983}, the \seq\ of \pp es
$\wt N_n=\sum_{i=1}^{k_n} \wt N_{nm}^{(i)}$ is relatively compact
and the subsequential limits  are infinitely divisible, possibly
null.  By virtue of
${\mathcal A}(a_n)$, $N_n\std
N$ for some infinitely divisible \pp\ $N$ \fif\ $\wt N_n\std N$.
\ere
\subsection{Weak \con\ of \pp es}
Now we are in the position  to formulate one of the  main results in
Davis and Hsing
\cite{davis:hsing:1995}. The result was proved
in the case $d=1$ but immediately translates to the case $d>1$; see
Davis and Mikosch
\cite{davis:mikosch:1998}.
\bth\label{thm:davhs}
Assume that the strictly stationary  $\bbr^d$-valued \seq\ $(X_t)$
satisfies
\begin{enumerate}
\item[\rm 1.]
the \regvar\ condition ${\bf (RV_\alpha)}$ for some $\alpha>0$, 
\item[\rm 2.]
the anti-clustering condition {\bf (AC)},
\item[\rm 3.]
the mixing condition
${\mathcal A}(a_n)$. 
\end{enumerate}
Then $N_n=\sum_{t=1}^n \vep_{a_n^{-1} X_t}\std N$ in the space of point \ms s on
$\ov\bbr_0^d$ equipped with the vague topology and the infinitely divisible
limiting \pp\ $N$ has \rep\
\beao
N=\sum_{i=1}^\infty\sum_{j=1}^\infty \vep_{\Gamma_i^{-1/\alpha} Q_{ij}}\,,
\eeao
$(\Gamma_i)$ is an increasing enumeration of a homogeneous Poisson
process
on $(0,\infty)$ with intensity $\gamma$, $(Q_{ij})_{j\ge 1}$,
$i=1,2,\ldots$, are iid \seq s of points $Q_{ij}$ \st\
$\max_{j\ge 1}|Q_{ij}|=1$ a.s. and $\gamma$ is the extremal index of
the \seq\ $(|X_t|)$. The \ds\ of $\sum_{j=1}^\infty \vep_{Q_{1j}}$ is
given in Theorem~2.7 of Davis and Hsing  \cite{davis:hsing:1995} in terms of the limit
\ms s  $\mu_k$, $k\ge 1$; see \eqref{eq:1a}.
\ethe
\bre\label{rem:bs}
Basrak and Segers \cite{basrak:segers:2009} gave an alternative proof
of this result. They also
showed that $N$ has
Laplace \fct al in terms of the spectral tail chain $(\Theta_t)$
given by
\beao
\E\ex ^{- \int f dN}= \exp\Big\{- \int_0^\infty
\E\Big(\ex^{-\sum_{i=1}^\infty f(y
  \Theta_i)}- \ex^{-\sum_{i=0}^\infty f(y \Theta_i)}\Big) 
 d(-y^{-\alpha}) \Big\}\,,\quad f\in\bbc_K^+\,.
\eeao
Moreover, they showed that the extremal index $\gamma$ is positive and
has \rep
\beam\label{eq:extremalindex}
\gamma=\E\Big(\sup_{t\ge 0}|\Theta_t|^\alpha-\sup_{t\ge 1}|\Theta_t|^\alpha\Big)\,.
\eeam
\ere

\section{A \ld s approach to limit theory for heavy-tailed time series}\label{sec:ld}
We provide
a general \ld\ result for \fct als acting on a \regvary\ \seq\ and
vanishing in some neighborhood of the origin. The latter property
means that only large values of the \seq\ make a contribution to the 
limiting quantities. 
\par
We consider a complex-valued  \fct\ $f$ on $(\bbr^d)^\bbn$ an denote its restrictions to 
 $(\bbr^d)^l$ by $f_l$, $f_0=0$. We will say that $f$ (or equivalently $(f_l)$) satisfies {\bf (C$_\vep$)}:
 \begin{itemize}
\item[{\bf (C$_\vep$)}] For $l\ge 1$, the a.e. continuous function $f_l$ is bounded uniformly  and  
$$
f_{l}(x_1,\ldots,x_{j_1-1},x_{j_1},\ldots,x_{j_2},x_{j_2+1},\ldots,x_l)
=f_{j_2-j_1+1}(x_{j_1},\ldots,x_{j_2})\,,\qquad1\le j_1\le j_2\le l\,,
$$
provided $|x_i|\le \vep$, $i=1,\ldots,j_1-1,j_2+1,\ldots,l$. 
\end{itemize}
We suppress the dependence on $\vep$ in the notation and we often 
write $f$ instead of $f_l$; it will be clear from the number of
arguments which $f_l$ we are dealing with.
\par
For a stationary  \seq\ $(X_t)$ we follow an approach which was
advocated by Jakubowski and Kobus \cite{jakubowski:kobus:1989} and Jakubowski \cite{jakubowski:1993,jakubowski:1997} in the
context of $\alpha$-stable limit theory for sums of infinite variance
\rv s and was exploited in Bartkiewicz et al.
\cite{bartkiewicz:jakubowski:mikosch:wintenberger:2011}, Balan and
Louhichi \cite{balan:louhichi:2009}, Mikosch and Wintenberger \cite{mikosch:wintenberger:2013,mikosch:wintenberger:2012}
for proving limit theory for \pp es, sums of \regvary\ \seq s with
infinite variance stable limit laws, \ld\ \pro ies and other results.
The main idea of this approach is to use suitable telescoping sums involving 
the differences $\E\big[ f(x^{-1}X_0,\ldots,x^{-1} X_k)-
f(x^{-1}X_1,\ldots,x^{-1} X_k)]$, $k\ge 1$.
\bth\label{th:antic}
Consider a strictly stationary $\bbr^d$-valued  \seq\ $(X_t)$
satisfying   
\begin{enumerate}
\item[\rm 1.]
${\bf (RV_\alpha)}$ for some $\alpha>0$,
\item[\rm 2.] 
the uniform anti-clustering condition
\beam\label{eq:antic}
\lim_{\kto}\limsup_{\nto} \sup_{x\in\Lambda_n} \P( M_{k,n}>x \delta\mid
|X_0|>x\delta )= 0\,,\quad \delta>0\,,
\eeam
for some \seq\ of Borel subsets $\Lambda_n\subset (0,\infty)$ \st\
$x_n=\inf \Lambda_n\to\infty$ and $n\,\P(|X_0|>x_n)\to 0$.
\end{enumerate}
Let $f$ be a  complex-valued \fct\ satisfying  {\bf
  (C$_\vep$)} for some $\vep>0$.
Then the following relation holds:
\beam\label{eq:exp}
\sup_{x\in\Lambda_n}\Big|
\frac{\E[ f(x^{-1}X_1,\ldots,x^{-1}X_{n})]}{n\P(|X_0|> x)}-
\int_{0}^\infty\E[f(y(\Theta_t)_{t\ge 0})-f(y(\Theta_t)_{t\ge
  1})]d(-y^{-\alpha})\Big|\to 0\,.\nonumber\\
\eeam
\ethe
The proof is given in Section~\ref{sec:proofofthm1}.
\bre\label{rem:ac+}
If we consider one-point sets $\Lambda_n$ the
anti-clustering condition \eqref{eq:antic} can be compared with
the anti-clustering condition {\bf (AC)} of Davis and Hsing
\cite{davis:hsing:1995}; see \eqref{eq:accond}. Indeed, if we choose 
$x_n=\delta a_n$ for $\delta>0$ and replace $M_{k,n}$ by $\wt M_{k,m_n}$
for some \seq\ $m_n\to\infty$, $m_n=o(n)$, then \eqref{eq:antic} is
related to {\bf (AC)}. However, there is one major difference: 
\eqref{eq:antic} does not involve maxima over sets of negative
integers. If the stronger condition {\bf (AC)} holds, Basrak and Segers
\cite{basrak:segers:2009} showed that the extremal index $\gamma_{|X|}$ of the \seq\ 
$(|X_t|)$ is positive. It is also the case under the less restrictive condition \eqref{eq:antic} because $Y_t\stp 0$ when $t\to\infty$ a.s. from the proof of Proposition 4.2  of \cite{basrak:segers:2009} and $\gamma_{|X|}=\P(|Y_t|\le 1,\, t>0)$. 
\ere

Recall the notion of a 
$k$-dependent stationary \seq\ $(X_t)$ for some integer $k\ge 0$, i.e.,
the $\sigma$-fields $\sigma(X_t,t\le 0)$ and $\sigma(X_t,t\ge k+1)$ 
are independent. In this case, Theorem~\ref{th:antic} simplifies.
\bth\label{th:kdep}
Consider a stationary $\bbr^d$-valued $k$-dependent  \seq\ $(X_t)$
satisfying  ${\bf (RV_\alpha)}$ for some $\alpha>0$ and a  
complex-valued \fct\ 
$f$ satisfying {\bf (C$_\vep$)} for some $\vep>0$. Let $(x_n)$ be a real-valued  \seq\ \st\ 
$n\,\P(|X_0|>x_n)\to0$. Then, as $\nto$,
\beam\label{eq:step1}
\sup_{x\ge x_n}\left|\frac{\E[
    f_{n}(x^{-1}X_1,\ldots,x^{-1}X_{n})]}{n\P(|X_0|> x)}-
  \la_k\int_{0}^\infty\E[f_{k+1}(y\wt\Theta_0,\ldots,y\wt\Theta_{k})]d(-y^{-\alpha})\right|\to 0\,.\nonumber\\
\eeam
with $\la_k=\P(\Theta_{-j}=0,j=1,\ldots,k)>0$ and 
\beam\label{eq:Theta}
\P\big( (\wt \Theta_{j})_{j=0,\ldots,k}\in\cdot\big)=
\P\big( ( \Theta_{j})_{j=0,\ldots,k}\in\cdot\mid \Theta_{-l}=0,
l=1,\ldots,k\big)\,.
\eeam
\ethe
The proof is given in Section~\ref{sec:proofofth2}.
\bre
The 
limiting expression in \eqref{eq:exp} 
 does in general not coincide with the ``naive'' limit by letting
 $k\to\infty$ in \eqref{eq:step1}, i.e.,  
$$
\P(\Theta_{-j}=0,j\ge 1)\int_{0}^\infty\E[f(y(\tilde \Theta_t)_{t\ge 0})]d(-y^{-\alpha})
$$
because $\P(\Theta_{-j}=0,j\ge 1)=0$ is possible. 
However, if $\P(|\Theta_{-l}|\le \delta, l\ge 1)>0$ for some
$\delta>0$
then we can define $(\wt \Theta_t^\delta)$ through the relation  \beao
\P\big( (\wt \Theta_{t}^{\delta})_{t\ge 0}\in\cdot\big)=
\P\big( ( \Theta_{t})_{t\ge 0}\in\cdot\mid |\Theta_{-l}|\le \delta,
l\ge 1\big)\,,\qquad\vep \ge \delta>0\,.
\eeao
In view of condition  {\bf (C$_\vep$)} and the proof of Theorem \ref{th:kdep}, we also have
\beao\lefteqn{
\int_{0}^\infty\E[f(y(\Theta_t)_{t\ge 0})-f(y(\Theta_t)_{t\ge
  1})]d(-y^{-\alpha})}\\
&=&\int_{\vep}^\infty\E[f(y(\Theta_t)_{t\ge 0})-f(y(\Theta_t)_{t\ge
  1})]d(-y^{-\alpha})\\
&=&\int_{\vep}^\infty\E f(y(\Theta_t)_{t\ge 0}) d(-y^{-\alpha})-
\int_{\vep}^\infty \E f(y(\Theta_t)_{t\ge
  1})]d(-y^{-\alpha})\,\\
  &=&\P(|\Theta_{-l}|\le \delta, l\ge 1)\int_{\vep}^\infty\E f(y(\tilde \Theta_t^\delta)_{t\ge 0}) d(-y^{-\alpha}),
\eeao
and both quantities on the \rhs\ in the third line are finite and involve only a finite number of $\Theta_t$ a.s. because $\Theta_t\stp 0$. 
\ere

\section{Applications}\label{subsec:appl}\setcounter{equation}{0}
In this section we will provide various applications of
Theorems~\ref{th:antic} and \ref{th:kdep}
to limit theorems
of \ld -type and weak \con\ results of various kinds.

\subsection{Limit theory for \pp es and partial sums}\label{subsec:limittheory} 
\subsubsection{Weak convergence of point processes}
We re-prove the \pp\ result of Davis and Hsing
\cite{davis:hsing:1995} on \pp\ \con\ given as Theorem~\ref{thm:davhs}
above, formulated in the language of Basrak 
and Segers \cite{basrak:segers:2009}; see Remark~\ref{rem:bs}. 
A careful analysis of \cite{basrak:segers:2009} shows that their proofs
use ideas which are close to those in the proof of
Theorem~\ref{th:antic}; see also Balan and 
Louhichi \cite{balan:louhichi:2009} who apply Jakubowski's ideas to
\pp\ \con\ of more general triangular arrays.
Recall the definition of the \pp es $N_n$ from \eqref{eq:nn}.
\bth\label{cor:ppgeneral}
Assume that the strictly stationary $\bbr^d$-valued   \seq\ $(X_t)$
satisfies  
\begin{enumerate}
\item[\rm 1.]
${\bf (RV_\alpha)}$ for some $\alpha>0$,
\item[\rm 2.] 
the mixing condition ${\mathcal A}(a_n)$,
\item[\rm 3.] 
the   anti-clustering
condition \eqref{eq:antic} for $\Lambda_n=\{a_n\}$. 
\end{enumerate}
Then $N_n\std N$,
where $(N_n)$ is defined in \eqref{eq:nn},
and $N$ has Laplace \fct al
\beao
\E\ex^{-\int _{\ov \bbr_0^d} g dN} = \exp\Big(-\int_{0}^\infty  \E (\ex^{-\sum_{j=1}^{\infty}
  g(y \Theta_j)}-\ex^{-\sum_{j=0}^{\infty}
  g(y \Theta_j)})\,d(-y^{-\alpha})\Big)\,,\quad  g\in \bbc_K^+\,.
\eeao
\ethe
\begin{proof}
We prove the convergence of the logarithms of the 
Laplace \fct als $\log \E \ex^{-\int
  gdN_n}\to  \log\E \ex^{-\int
  gdN}$ for $g\in\bbc_K^+$. 
We assume that $g(x)=0$ for
$|x|\le \vep$ for some $\vep>0$. In view of the mixing condition
we have for some $m=m_n\to\infty$ and $k_n=[n/m]\to\infty$,
\beao
\log \E\ex ^{- \int gdN}=\log \E\ex^{-\sum_{t=1}^ng(a_n^{-1}X_t)}\sim k_n  \log \E\ex^{-\sum_{t=1}^{m_n}g(a_n^{-1}X_t)}\,.
\eeao
By a Taylor expansion, since $\E\sum_{t=1}^{m_n}g(a_n^{-1}X_t)\le C
m_n\P(|X|>\vep a_n)\to 0$,
\beao
- k_n  \log\E \exp(-\sum_{t=1}^{m_n}g(a_n^{-1}X_t))\sim k_n \E(f_{m_n}(a_n^{-1}X_1,\ldots,a_n^{-1}X_{m_n})),
\eeao
where 
\beao
f_{l}(x_{1},\ldots,x_{l})= 1-\exp\Big(-\sum_{t=1}^{l}g(x_t)\Big)\,,\quad 1\le l\,.
\eeao
Notice that $(f_l)$ satisfies {\bf (C$_\vep$)}.
Now an  application of Theorem \ref{th:antic} with $\Lambda=\{a_n\}$
and $n$ replaced by $m_n$ yields 
\beao
k_n\,  \E(f(a_n^{-1}X_1,\ldots,a_n^{-1}X_{m_n}))\to \int_0^\infty 
\big[\E [f(y(\Theta_{i})_{i\ge 0})- f(y(\Theta_{i})_{i\ge 1})\big]\,d(-y^{-\alpha})\,.
\eeao
The limit is the desired logarithm of the Laplace \fct al $N$.
Combining the arguments, we proved the corollary.
\end{proof}
\bre
For a $k$-dependent \regvary\ \seq\ $(X_t)$, the mixing and
anti-clustering
conditions of Theorem~\ref{cor:ppgeneral} are trivialy satisfied. Moreover,
we conclude from Theorem~\ref{th:kdep} that $N$ has Laplace \fct al
\beao
E\ex^{-g dN} = \exp\Big(-\la_k\,\int_{0}^\infty  \E (1-\ex^{-\sum_{j=0}^{k}
  g(y \wt \Theta_j)})\,d(-y^{-\alpha})\Big)\,,\quad  g\in \bbc_K^+\,.
\eeao
Calculation shows that
the infinitely divisible
limiting \pp\ $N$ has \rep\
\beao
N=\sum_{i=1}^\infty\sum_{j=0}^{k} \vep_{\Gamma_i^{-1/\alpha} \wt \Theta_{ij}}\,,
\eeao
where $(\Gamma_i)$ is an increasing enumeration of a homogeneous Poisson
process on $(0,\infty)$ with intensity
$\la_k$, independent of an iid \seq\ 
$( \wt \Theta_{ij})_{0\le j\le k}$, $i=1,2,\ldots$, with 
generic element $( \wt \Theta_{j})_{0\le j\le k}$. 
\ere
\subsubsection{The $\alpha$-stable central limit theorem} 
In this section we consider the truncated \rv s 
\beao
\ov X_t= X_t \1_{|X_t|\le \vep a_n}\,,\quad \underline X_t=X_t-\ov
X_t\,,\quad t\in\bbz\,,  
\eeao
and the corresponding partial sums $\ov S_n$ and $\underline S_n$,
where we suppress the dependence on $\vep>0$ and $n$ in the notation.
\bth\label{thm:alphaclt}
Consider  a stationary $\bbr^d$-valued \seq\ $(X_t)$ satisfying
 \begin{enumerate}
\item[\rm 1.] 
 ${\bf (RV_\alpha)}$ for some $\alpha\in (0,2)\backslash\{1\}$,
\item[\rm 2.] 
the mixing condition
\beam\label{eq:mix}
\E \ex^{is' \underline S_n/a_n}- \Big(\E \ex^{is' \underline
  S_m/a_n}\Big)^{k_n}\to 0\,,\quad s\in\bbr^d\,,\quad \nto\,,
\eeam
\item[\rm 3.] 
the anti-clustering condition \eqref{eq:antic} for $\Lambda_n=\{a_n\}$,
\item[\rm 4.] 
for $\alpha\in (1,2)$, in addition, $\E X=0$, the vanishing-small-values condition
\beam\label{eq:neg}
\lim_{\vep\downarrow 0}\limsup_{\nto} \P(a_n^{-1} |\ov
S_n -\E \ov S_n|>\delta)=0\,,\quad \delta>0\,,
\eeam 
and $\sum_{j=1}^\infty \E |\Theta_j|<\infty$.
\end{enumerate} 
Then $a_n^{-1} S_n\std \xi_\alpha$, where
the limit is an $\alpha$-stable \rv\ with log-\chf 
\beam\label{eq:stable}
\int_0^\infty \E\Big( \ex^{i y s' \sum_{j=0}^\infty \Theta_j }-
\ex^{i y s' \sum_{j=1}^\infty \Theta_j}-iys' A\Big)\,d (-y^{-\alpha})\,, \qquad s\in\bbr^d\,,
\eeam
where $A=0$ for $\alpha\in (0,1)$ and $A=\Theta_0$ for $\alpha\in
(1,2)$.
\ethe
The proof is given in Section~\ref{sec:proofofth}.
\bre
Condition~\eqref{eq:neg} for $\alpha\in (1,2)$ is standard in central
limit theory; see the discussions in 
Davis and Hsing \cite{davis:hsing:1995},
Bartkiewicz et
al. \cite{bartkiewicz:jakubowski:mikosch:wintenberger:2011},
Basrak et al. \cite{basrak:krizmanic:segers:2012}. Sufficient
conditions are $k$-dependence of $(X_t)$ and conditional independence.
For concrete models such as stochastic volatility models, GARCH and
certain Markov chains, see the references above and 
\cite{mikosch:wintenberger:2013,mikosch:wintenberger:2012}.
For similar characterizations of the $\alpha$-stable limiting laws as in \eqref{eq:stable}, see Mirek \cite{mirek:2011}. It coincides with the limit law given in \cite{bartkiewicz:jakubowski:mikosch:wintenberger:2011} as shown by the computations of Louhichi and Rio \cite{louhichi:rio:2011}.

In the $k$-dependent case, Jakubowski and Kobus \cite{jakubowski:kobus:1989} and Kobus \cite{kobus:1995} got related $\alpha$-stable limit theory under the assumption that $(X_0,\ldots,X_k)$ is \regvary\ with index $\alpha$. In view of Proposition \ref{lem:1}, the latter condition is equivalent to condition {\bf (RV$_\alpha$)}. Extensions of the $\alpha$-stable \clt\ to the stationary case were considered in Jakubowski \cite{jakubowski:1993,jakubowski:1997}.
 \ere

\subsection{Large deviations for suprema of a random walk and ruin bounds}\label{sebsec:supremaandruin}

\subsubsection{Large deviations for the  supremum of a random walk}\label{sec:ldsup}
In this section we derive a result for the suprema of a univariate 
random walk $(S_n)$. 
We write  
for any $x,\vep>0$,
\beao
\ov X_t = X_t\1_{|X_t|\le \vep x}\,,\quad \underline X_t= X_t- \ov
X_t\,,\quad t\in\bbz\,,  
\eeao
and 
\beao
S_0=0\,,\quad  \ov S_t= \sum_{i=1}^t \ov X_i\,,\quad\underline
S_t=S_t-\ov S_t\,,\quad t\in \bbz\,.
\eeao
Here we suppress the dependence of these quantities on $x,\vep$ in the
notation.
\bth\label{thm:suprema}
Consider  a stationary $\bbr$-valued \seq\ $(X_t)$ satisfying
the following conditions
\begin{enumerate}
\item[\rm 1.] 
 ${\bf (RV_\alpha)}$ for some $\alpha>0$,
\item[\rm 2.]
the anti-clustering condition \eqref{eq:antic}.
\end{enumerate}
If $\alpha>1$ we also assume
\begin{enumerate}
\item[\rm 3.]
the vanishing-small-values condition
\beam\label{eq:nega}
\lim_{\vep\downarrow 0}\limsup_{\nto} \sup_{x\in \Lambda_n}\dfrac{\P(x^{-1} \sup_{t\le n}|\ov
S_t  |>\delta)}{n\,\P(|X|>x)}=0\,,\quad \delta>0\,,
\eeam 
\item[\rm 4.]
$\E \Big(\sum_{i=1}^\infty
|\Theta_i|\Big)^{\alpha-1}<\infty$.
\end{enumerate} 
Then 
\beam\label{eq:supsn}
\qquad \sup_{x\in \Lambda_n}\Big|\dfrac{\P\Big(\sup_{i\le n}
  S_i>x)\Big)}{n\,\P(|X|>x)} - \E\Big[\Big( \Theta_0+\sup_{t\ge
  1}\sum_{i=1}^t\Theta_i\Big)_+^\alpha - \Big(\sup_{t\ge
  1}\sum_{i=1}^t\Theta_i\Big)_+^\alpha\Big]
\Big|\to 0\,,\quad \nto\,.
\eeam
\ethe
The proof is given in Section~\ref{prof}.

\bre
We observe that the limit in \eqref{eq:supsn} can be $0$, for example for the large deviations of the telescoping sum $S_n$ with $X_t=Y_t-Y_{t-1}$ with $Y_t\ge0$ iid \regvary. Then $\Theta_0=-\Theta_1=1$, $\Theta_t=0$ for $t\ge2$ and the limiting constant is zero.
\ere
Condition \eqref{eq:nega} can often be verified by using maximal
inequalities for sums, for example
in the case of regenerative \MC s or conditionally 
independent random variables; see for example \cite{mikosch:wintenberger:2012}.
For a $k$-dependent sequence one can verify this condition as well (see 
Lemma \ref{lem:8}), resulting in the following corollary.
\bco\label{th:sup sum}
Assume that $(X_t)$ is a $k$-dependent univariate strictly 
stationary \seq\ which is also  \regvary\
with index $\alpha>0$. In addition, we assume the following
conditions:
\begin{enumerate}\item[\rm 1.]
$\E X=0$ if $\E|X|<\infty$.
\item[\rm 2.] 
If $\alpha=1$ and $\E |X|=\infty$, we have
\beam\label{eq:auxo}
\limsup_{\nto}\sup_{x\ge x_n} n x^{-1} |\E\ov X|=0\,.\eeam
\end{enumerate}
Let $(a_n)$ be any \seq\ \st\ $n\,\P(|X|>a_n)\to 1$ as $\nto$.
Then 
\beao
\sup_{x\ge x_n}\Big|\frac{ \P\Big(\sup_{t\le n}   S_t >x\Big)}{n\P(|X|>x)}-\lambda_k\E\Big( \sup_{t\le k}
\sum_{i=0}^t\wt\Theta_i \Big)_+^\alpha\Big|\to0,
\eeao
where $x_n\to\infty$ is any \seq\ \st\ $x_n/a_n\to\infty$ if $\alpha<2$, $x_n/n^{0.5+\delta}\to\infty$
for some $\delta>0$ 
if $\alpha=2$ and $\E X^2=\infty$, and  
$x_n\ge C\sqrt{n\log n}$ for sufficiently large $C>0$ if $\E
X^2<\infty$.
\eco
The proof is given in Section~\ref{proofofcorol}.
\par
\bre
The proof of Corollary~\ref{th:sup sum} immediately extends to certain 
subadditive
\fct als acting on
the random walk $(S_n)$ which are more general than suprema. Indeed,
let $g_l$ be a \seq\ of real-valued \fct s on $\bbr^l$, $l\ge 1$.
Assume that, for any $l\ge 1$,
\begin{itemize}
\item $g_l$ is
{\em continuous and positively homogeneous,} i.e., $g_l(c\bfx)=cg_l(\bfx)$ 
for any $\bfx\in\bbr^l$ and $c>0$,
\item 
{\em subadditive,} i.e.,
$
g_l(\bfx+\bfy)\le g_l(\bfx)+
g_l(\bfy),\quad \bfx,\bfy\in\bbr^l\,,
$
\item
a {\em domination property} holds: for $s_t=x_1+\cdots
+x_t$, $\bfs_l=(s_1,\ldots,s_l)$, there exists a constant
$C>0$ not depending on $l$ \st
$
|g_l(\bfs_l)|\le C \sup_{ 1\le i\le l}|s_i|\,. 
$ 
\end{itemize}
Finally, write $\underline s_t=\sum_{i=1}^t x_i \1_{|x_i|>\vep }$
and assume that $(\1_{g_l(\underline s_1,\ldots,\underline s_l)>1})$
satisfies {\bf (C$_\vep$)}.
Then, under the conditions and with the notation of Corollary~\ref{th:sup sum}, the following result
holds:
\beao
\sup_{x\ge x_n}\Big|\frac{ \P\Big(g_n\big((S_t)_{t=1,\ldots,n}\big)
  >x\Big)}{n\P(|X|>x)}-\lambda_k\E\Big[\Big( g_{k+1}\big(\big(
\sum_{i=0}^t\wt\Theta_i\big)_{t=0,\ldots,k}\big) \Big)_+^\alpha
\Big]\Big|\to0\,.
\eeao
For example, with $g_l(s_1,\ldots,s_l)=\sup_{t\le l} |s_t|$, $l\ge 1$, we obtain
\beao
\sup_{x\ge x_n}\Big|\frac{ \P\Big(\sup_{t\le n}   |S_t| >x\Big)}{n\P(|X|>x)}-\lambda_k\E \sup_{t\le k}\Big|
\sum_{i=0}^t\wt\Theta_i \Big|^\alpha\Big|\to0\,.
\eeao
With  $g_l(s_1,\ldots,s_l)=s_l$ we 
get a \ld\ result for sums: 
\beao
\sup_{x\ge x_n}\Big|\frac{ \P\Big(S_n >x\Big)}{n\P(|X|>x)}-\lambda_k\E\Big( 
\sum_{i=0}^{k}\wt \Theta_i\Big)_+^\alpha\Big|\to0\,.
\eeao
Other \fct als of this kind are given by
$g_l(s_1,\ldots,s_l)=\max_{i=1,\ldots,l}(s_i - s_l)_+$ and 
$g_l(s_1,\ldots,s_l)=\max_{i=1,\ldots,l}s_i - 
\min_{i=1,\ldots,l}s_i$, $g_l(s_1,\ldots,s_l)= \max_{i,j=1,\ldots,l} (s_i-s_j)$.

For \regvary\ moving averages with index $\alpha<2$, Basrak and Krizmani\`c \cite{basrak:krizmanic:2014} studied the $M_2$-\fct al convergence of the partial sums when $\sum_{i=0}^{k}\wt\Theta_i$ coincides with $\sup_{0\le t\le k}\sum_{i=0}^t\wt\Theta_i$. They derived the limiting law of the supremum of the partial sums; it is the supremum of the $\alpha$-stable L\'evy process $(\xi_t)_{t\in [0,1]}$ where $\xi_1$ has the same distribution as the limiting law of the partial sums. In view of the results above, a similar phenomenon can be observed under the condition
\beam\label{eq:m2}
\E\Big( 
\sum_{i=0}^{k}\wt \Theta_i\Big)_+^\alpha=\E\Big( \sup_{0\le t\le k}
\sum_{i=0}^{t}\wt \Theta_i\Big)_+^\alpha
\eeam
for the large deviations of sums and their suprema:
$$
\sup_{x\ge x_n}\Big|\frac{\P(S_n >x)}{n\P(|X|>x)}- \frac{\P(\sup_{1\le t\le n}S_t >x)}{n\P(|X|>x)}\Big|\to 0\,.
$$
However, in the cases where \eqref{eq:m2} does not hold, the functional \clt\ cannot hold for any topology for which the supremum is a continuous function.
\ere
\subsubsection{Ruin probabilities}\label{sec:ruina}
In this section we assume that $(X_t)$ is a univariate strictly
stationary \seq\ 
which is also \regvary\ with index $\alpha>1$. The
latter condition ensures that $\E|X|<\infty$. We will also assume that
$\E X=0$. In what follows, we will study the \asy\ behavior of the
tail \pro y
$\P(\sup_{t\ge0}(S_t-\rho t)>x)$ for $\rho>0$ as $\xto$. We will refer
to this \pro y as {\em ruin \pro y} since similar expressions appear
in the context of non-life \ins\ mathematics; see Asmussen and Albrecher
\cite{asmussen:albrecher:2010} and Embrechts et al.
\cite{embrechts:kluppelberg:mikosch:1997}, Chapter 1.
We use the notation of Section~\ref{sec:ldsup}.
\bth\label{thm:ruina}
Assume that 
$(X_t)$ is a univariate strictly
stationary  \seq\ 
which is also \regvary\ with index $\alpha>1$, has mean zero and satisfies the conditions  of Theorem \ref{thm:suprema} with $\Lambda_n=[C_1n,C_2n]$ for any possible choice of positive constants $C_1<C_2$. Then 
we have for any $\rho>0$,
\beam\label{eq:ruin}
\frac{\P(\sup_{t\ge0}(S_t-\rho t)>x)}{x\,\P(|X|>x)}\sim \frac
{\E\Big[\Big( \sup_{t\ge 0}\sum_{i=0}^t  \Theta_i\Big)_+^\alpha-\Big( \sup_{t\ge1}\sum_{i=1}^t  \Theta_i\Big)_+^\alpha\Big]}{(\alpha-1)\rho}\,,\quad \xto \,.\nonumber\\
\eeam
\ethe
The proof is given in Section~\ref{proofofruin}.
\bre\label{rem:pos}
Notice that the \rhs\ of \eqref{eq:ruin} is of the form
\beao
\frac
{\E [ ( 1+\sum_{i=1}^\infty  \Theta_i )^\alpha- ( \sum_{i=1}^\infty  \Theta_i)^\alpha ]}{(\alpha-1)\rho}
\eeao
provided that the $\Theta_t$, $t\ge 0$, are non-negative.
\ere

\bco\label{thm:ruin}
Assume that 
$(X_t)$ is a univariate strictly
stationary $k$-dependent \seq\ 
which is also \regvary\ with index $\alpha>1$ and has mean zero. Then 
we have for any $\rho>0$,
\beao
\frac{\P(\sup_{t\ge0}(S_t-\rho t)>x)}{x\,\P(|X|>x)}\sim \frac
{\lambda_k}{(\alpha-1)\rho}
\E\Big( \sup_{t\le
    k}\sum_{i=0}^t\wt \Theta_i\Big)_+^\alpha\,,\quad \xto \,.
\eeao
\eco
The proof is given in Section~\ref{proofofruina}.
\bexam\label{exam:goldie}\em
Consider the \sre\ $X_t=A_tX_{t-1}+B_t$, $t\in\Z$, where $(A_t,B_t)$,
$t\in\Z$, constitute an $\R^2$-valued iid sequence. We assume that $(X_t)$ constitutes a strictly stationary \MC .\\

{\it The Goldie case:}  We assume
the conditions of Goldie \cite{goldie:1991} are satisfied, ensuring
that $X_0$ is \regvary\ with index   $\alpha>0$. In particular, we
have $A\ge 0$ a.s.,
$\E[A^\alpha]=1$ for some positive $\alpha$ and we also need some
further conditions on the distribution of the sequence $(A_t,B_t)$ 
to ensure that one has a Nummelin regeneration scheme. Then $(X_t)$ is
\regvary\ of order $\alpha$,
$\Theta_0=1$ and  $\Theta_t=\Pi_t=A_1\cdots A_t$, $t\ge 1$; see Basrak
and Segers \cite{basrak:segers:2009}.  Proceeding 
as in the proof of Theorem 7.2 in \cite{mikosch:wintenberger:2012} and
using the drift condition\footnote{We say that a function $X_t=g(\Phi_t)$ of a Markov chain $(\Phi_t)$ satisfies the {\bf (DC$_p$)} condition if $\E[|g(\Phi_1)|^p\mid \Phi_0=x]\le \beta|g(x)|^p+b\1_{A}(x)$ where $0<\beta<1$ and $A$ is a small set; see \cite{meyn:tweedie:1993} for details.} {\bf (DC$_p$)} with $p<\alpha$, 
one can show the anti-clustering condition \eqref{eq:antic}. 
In the proof of Theorem 4.6 in \cite{mikosch:wintenberger:2013} 
we showed that the vanishing-small-values condition \eqref{eq:vsv} without the supremum
is satisfied under {\bf (DC$_p$)} for  $p<\alpha$. 
An inspection of the proof also shows that one may restrict oneself 
to the absolute values $|\ov X_t|$, implying the vanishing-small-values condition for the supremum as well. From Theorem \ref{thm:ruina} and Remark \ref{rem:pos} we conclude that 
\beam\label{eq:tailgldie}
\frac{\P(\sup_{t\ge0}(S_t-(\rho+ \E X) t)>x)}{x\,\P(|X|>x)}\sim \frac
{\E\Big[\Big( 1+\sum_{i=1}^\infty  \Pi_i\Big)^\alpha-\Big(  \sum_{i=1}^\infty \Pi_i\Big)^\alpha\Big]}{(\alpha-1)\rho}\,,\quad \xto \,.\nonumber\\
\eeam
This result recovers Theorem 4.1 in Buraczewski et
al. \cite{buraczewski:damek:mikosch:zienkiewicz:2013} in  the case of
a Nummelin regeneration scheme. The method of proof in
\cite{buraczewski:damek:mikosch:zienkiewicz:2013} is completely
different  from ours. 
\par
We also mention that  
$Y_t\eqd 1+\sum_{i=1}^\infty  \Pi_i$, where $(Y_t)$ is the strictly
stationary causal solution to the \seq\ $Y_t=A_t Y_{t-1}+1$,
$t\in\bbz$, which is \regvary\ with index $\alpha$. Then the constant on the
\rhs\ of 
\eqref{eq:tailgldie} can be written as 
\beao
\frac
{\E\Big[Y_0^\alpha-(Y_0-1)^\alpha\Big]}{(\alpha-1)\rho}\,.
\eeao

{\it The Grey case:} We assume now that the conditions of Grey
\cite{grey:1994} are satisfied. This means that $A$ may assume real values,
$\E|A|^\alpha<1$ and $B$ \regvary\ with index $\alpha$ for some $\alpha>0$. Then the unique strictly stationary solution to the \sre\ exists and is \regvary\ with index $\alpha$ and $\Theta_t/\Theta_0=\Pi_t$ as above. Following Segers \cite{segers:2007}, we have in this case $\P(\Theta_{-j}=0,\,j\ge 1)=\P(\Theta_{-1}=0)=1-\E|\Theta_1|^\alpha>0$ because $|\Theta_1|=|\Theta_0||A_1|=|A_1|$. Thus, one can turn to the simpler alternative expression of the \rhs\ of \eqref{eq:ruin}
$$
(1-\E|A|^\alpha)\frac
{\E\Big[\sup_{t\ge 0}\Big( \sum_{i=0}^t  \wt \Theta_i\Big)_+^\alpha\Big]}{(\alpha-1)\rho},
$$
where $\wt \Theta_t/\wt \Theta_0=\Pi_t$, $\P(\wt \Theta_0=1)=\lim_{\xto}\P(B>x)/\P(|B|>x)=p$ and $\P(\wt \Theta_0=-1)=\lim_{\xto}\P(B\le -x)/\P(|B|>x)=q$. Then, we obtain 
$$
(1-\E|A|^\alpha)\frac
{\E\Big[p\sup_{t\ge 0}\Big(1+ \sum_{i=1}^t  \Pi_i\Big)_+^\alpha+q\sup_{t\ge 0}\Big(1+ \sum_{i=1}^t  \Pi_i\Big)_-^\alpha\Big]}{(\alpha-1)\rho}\,.
$$
We recover the result of Konstantinides and Mikosch \cite{konstantinides:mikosch:2005} for $A\ge 0$ a.s., $p=1$ and the one of Mikosch and Samorodnitsky \cite{mikosch:samorodnitsky:2000} in the AR(1) case when $A=\phi$ for some $|\phi|<1$. If $A\ge 0$ then the constant turns into 
$$
(1-\E A^\alpha)\frac
{p\,\E\Big[\Big(1+ \sum_{i=1}^t  \Pi_i\Big)^\alpha\Big]}{(\alpha-1)\rho}\,.
$$
Ruin bounds in the case of general linear processes $X_t=\sum_j\psi_jZ_{t-j}$ for iid \regvary\ $(Z_t)$ can be derived in a similar fashion using the computations of the spectral tail process given in Meinguet and Segers \cite{meinguet:segers:2010} recovering the results in \cite{mikosch:samorodnitsky:2000}.
\eexam
\bexam\label{exam:garch}\em
Consider the GARCH(1,1) model $X_t=\sigma_tZ_t$, $t\in\bbz$; see
Bollerslev \cite{bollerslev:1986}. Here $(Z_t)$ is an iid mean zero
and unit variance \seq\ of \rv s and $(\sigma_t^2)$ satisfies the \sre\
$\sigma_t^2=\alpha_0+(\alpha_1Z_{t-1}^2+\beta_1)\sigma_{t-1}^2$,
$t\in\Z$, 
where $\alpha_0$, $\alpha_1$ and $\beta_1$ are positive constants
chosen \st\ $(\sigma_t^2)$ is strictly stationary.  Moreover, we
assume that the above \sre\ for $(\sigma_t^2)$ with  $B_t=\alpha_0$
and $A_t=\alpha_1Z_{t-1}^2+\beta_1$ satisfies the Goldie conditions of 
Example~\ref{exam:goldie}, ensuring that $(\sigma_t^2)$ is \regvary\
with index $\alpha/2$. Then an application of Breiman's multivariate
results (see Basrak et al. \cite{basrak:davis:mikosch:2002}) implies that $(X_t)$
is \regvary\ with index $\alpha$. As before, we write
$\Pi_t=A_1\cdots A_t$.
Following \cite{mikosch:wintenberger:2012}, Section 5.4, 
we observe that as $\xto$,
\beao
\dfrac{\P(|(X_0,\ldots,X_t)-\sigma_0
(Z_0,\Pi_1^{0.5}
Z_1,\ldots,\Pi_t^{0.5}Z_t)|>x)}{\P(\sigma>x)} =o(1)\,.
\eeao
An application of the multivariate Breiman result yields 
$$
\frac{\P(x^{-1}\sigma_0
(Z_0,
Z_1\Pi_1^{0.5},\ldots,Z_t\Pi_t^{0.5})\in \cdot)}{\P(|X|>x)}\stv \frac{1}{\E|Z_0|^\alpha}\int_0^\infty \P(y(Z_0,Z_1\Pi_i^{0.5},\ldots,Z_t\Pi_t^{0.5})\in\cdot)\,d(-y^{-\alpha}).
$$
Then
\beao
\lefteqn{\P(x^{-1}(X_0,\ldots,X_t)\in \cdot\mid |X_0|>x)}\\
&\stw &\frac{1}{\E|Z_0|^\alpha}\int_0^\infty
\P(y(Z_0,Z_1\Pi_i^{0.5},\ldots,Z_t\Pi_t^{0.5})\in\cdot\,, y |Z_0|>1)\,d(-y^{-\alpha})\\
&=&\frac{1}{\E|Z_0|^\alpha}\E\Big[|Z_0|^\alpha \1_{|Y_0|(Z_0,Z_1\Pi_i^{0.5},\ldots,Z_t\Pi_t^{0.5})\in\,|Z_0|\,\cdot} \Big]\,,
\eeao
where $|Y_0|$ is Pareto distributed with index $\alpha$ and 
independent of $(Z_t)$.
By direct calculation, we obtain 
\beao
\lefteqn{
 \frac
{\E\Big[\Big( \sup_{t\ge 0}\sum_{i=0}^t  \Theta_i\Big)_+^\alpha-\Big( \sup_{t\ge1}\sum_{i=1}^t  \Theta_i\Big)_+^\alpha\Big]}{(\alpha-1)\rho}}\\
&=&\frac{\E\Big[
  \sup_{t\ge0}(Z_0+\sum_{i=1}^tZ_i\Pi_i^{0.5})_+^\alpha-\sup_{t\ge1}(\sum_{i=1}^tZ_i\Pi_i^{0.5})_+^\alpha\Big]}{\E|Z_0|^\alpha(\alpha-1)\rho}\,.
\eeao
Thus we derived the scaling constant for the ruin \pro y in \eqref{eq:ruin}
in the case of a \garch\ process. We also mention that the other
conditions of Theorem~\ref{thm:ruina} are satisfied. Indeed, the drift
{\bf (DC$_p$)} for $p<\alpha$ is satisfied, implying the
anti-clustering and vanishing-small-values conditions as in the case of
Example~\ref{exam:goldie}; see \cite{mikosch:wintenberger:2012} for details.
\eexam
\bexam\label{exam:sv}\em 
Consider a stochastic volatility model $X_t=\sigma_tZ_t$, $t\in\Z$, where $(\sigma_t)$ is a strictly stationary sequence with lognormal marginals independent of an iid \regvary\ sequence $(Z_t)$. Then $(X_t)$ is \regvary\ with the same index and it is not difficult to see that $\Theta_t=0$ for $t\neq0$. Now an application of Theorem \ref{thm:ruina} yields the same ruin bound as in the iid case. This result supports the general theory of such models whose extremal behaviour mimics the one of an iid \regvary\ sequence.
\eexam

\subsubsection{Large deviations for multivariate sums on half-spaces}
The same techniques as in the previous section can be used to prove
the following \ld\ result for multivariate sums.
\bth\label{thm:large}
Consider  a stationary $\bbr^d$-valued \seq\ $(X_t)$ satisfying
the following conditions
\begin{enumerate}
\item[\rm 1.] 
 ${\bf (RV_\alpha)}$ for some $\alpha>0$,
\item[\rm 2.]
the anti-clustering condition \eqref{eq:antic}.
\end{enumerate}
If $\alpha>1$ we also assume
\begin{enumerate}
\item[\rm 3.]
the vanishing-small-values condition
\beam\label{eq:vsv}
\lim_{\vep\downarrow 0}\limsup_{\nto} \sup_{x\in
  \Lambda_n}\dfrac{ \P (x^{-1} |\ov
S_n  |>\delta)}{n\,\P(|X|>x)}=0\,,\quad \delta>0\,,
\eeam
\item[\rm 4.]
$\E \Big(\sum_{i=1}^\infty
|\Theta_i|\Big)^{\alpha-1}<\infty$.
\end{enumerate} 
Then for every $\theta\in \bbs^{d-1}$,
\beao
\sup_{x\in \Lambda_n}\Big|\dfrac{ \P\Big(
  \theta' S_n>x)\Big)}{n\,\P(|X|>x)} - \E\Big[\Big( \theta'\sum_{i=0}^\infty\Theta_i\Big)_+^\alpha - \Big(\theta'\sum_{i=1}^\infty\Theta_i\Big)_+^\alpha\Big]
\Big|\to 0\,,\quad \nto\,.
\eeao
In addition, 
if $\alpha\not\in \bbn$ or $X$ is symmetric, then there exists a
unique Radon \ms\ $\mu_\alpha$ on $\ov\bbr_0^d$ \st\
$\mu_\alpha(t\cdot)=t^{-\alpha} \mu_\alpha(\cdot)$, $t>0$, and for any sequence $(x_n)$ such that $x_n\in\Lambda_n$, $n\ge 1$,
\beao
\dfrac{ \P\ (
  x_n^{-1}S_n\in \cdot  )}{n\,\P(|X|>x_n)}\stv
\mu_\alpha(\cdot)\,,\quad \nto\,.
\eeao
Moreover, $\mu_\alpha$ is determined by its values on the subsets
 $\{y\in\bbr^d: \theta'y>1\}$ for any $\theta\in\bbs^{d-1}$.
\ethe
The proof  of the last part follows by the same arguments as for Theorem 4.3 in 
Mikosch and Wintenberger \cite{mikosch:wintenberger:2012}.
\subsection{Large deviations  for cluster functionals}\label{subsec:clfctal}
Following  Yun \cite{yun:2000}  and Segers \cite{segers:2005}, we call a \seq\ of non-negative \fct s
$(c_l)$ on $\bbr^l$ a
{\em cluster \fct al} if it satisfies  {\bf (C$_0$)}.
As for the restrictions
$f_l$, we will often suppress the dependence of $c_l$ on the index $l$; it will be
clear from the context.
\par
Simple examples of cluster functionals are 
\beao
c_l(x_1,\ldots,x_l)&=&\sum_{i=1}^l\phi(x_i)\,,\quad l\ge 1\,,\eeao
with $\phi~: \R^d\to\R^+$ satisfying $\phi(0)=0$ and, for $d=1$,
\beao
c_l(x_1,\ldots,x_l)&=&\sum_{i=1}^l(x_i-z)_+\,,\quad l\ge 1\,,\\
c_l(x_1,\ldots,x_l)&=&\max_{1\le i\le l}(x_i-z)_+\,,\quad  l\ge 1\,,
\eeao
for some $z\ge 0$; see \cite{drees:rootzen:2010,yun:2000,segers:2005} for further examples.
\bco\label{co:clf}
Assume that the strictly stationary $\bbr$-valued   \seq\ $(X_t)$
satisfies  ${\bf (RV_\alpha)}$, the uniform anti-clustering
condition \eqref{eq:antic}  and that $$f_l(x_1,\ldots,x_l)=\1_{c_l((x_1-1)_+,\ldots,(x_l-1)_+)>1}$$ satisfies
{\bf (C$_1$)}. 
Then 
\beam\label{eq:lo0}
\lefteqn{\sup_{x\in \Lambda_n}\Big|\frac{\P(c((x^{-1}X_i-1)_+)_{1\le i\le
      n})>1)}{n\P(|X|>x)}}\nonumber\\&&-[\P(c(((|Y_0|\Theta_t-1)_+)_{t\ge0})>1)-\P(c(((|Y_0|\Theta_t-1)_+)_{t\ge1})>1)\Big|\to0\,,\quad \nto\,.
\eeam
\eco
\begin{proof}
We apply Theorem \ref{th:kdep} to $(f_l)$ satisfying
{\bf (C$_1$)} and we obtain the uniform limit for $x\in \Lambda_n$ as $\nto$:
\beao
\int_0^\infty\big[\P(c(((y \Theta_i-1)_+)_{i\ge 0})>1)
-\P(c(((y \Theta_i-1)_+)_{i\ge 1})>1)\big]
\,d(-y^{-\alpha}).
\eeao
By assumption,  $|\Theta_0|=1$ and therefore 
we can restrict the area of integration to $[1,\infty)$, recovering  
the desired limit. 
\end{proof}
\bre
In the $k$-dependent case, the anti-clustering condition is trivially
satisfied. In view of Theorem~\ref{th:kdep} relation \eqref{eq:lo0} then turns into 
\beao
\sup_{x\ge x_n}\Big|\frac{\P\Big(c_n\big(([x^{-1}X_i-1]_+)_{1\le i\le n}\big)>1\Big)}{n\,\P(|X|>x)}-\la_k\,\P\Big(c_{k+1}\big(([\,|Y_0|\wt\Theta_i-1]_+)_{0\le i\le k}\big)>1\Big)\Big|\to0\,.
\eeao
\ere
\bre
Corollary~\ref{co:clf} is formulated for univariate \seq s
$(X_t)$. However, one can generalize this result in various ways; see for example Drees and Rootz\'en \cite{drees:rootzen:2010}. For
example, let $(X_t)$ be an $\bbr^d$-valued \seq\ satisfying the
conditions of Theorem~\ref{th:kdep} and $A\subset \ov \bbr^d_0$ be a
Borel set whose distance to the origin is $\ge \vep>0$. Moreover, let
$g: \bbr^d\to \bbr$ be a measurable \fct .  If $g$ is a.e. continuous and $A$ is a continuity set with respect to Lebesgue measure then the same argument as for 
Corollary~\ref{co:clf} now yields
\beao\lefteqn{
\sup_{x\in \Lambda_n}\Big|\frac{\P(c_n( g(x^{-1}X_i) \1_{x^{-1} X_i\in A
    })_{1\le i\le n})>1)}{n\P(|X|>x)}}\\&&-\big[\P(c((g(|Y_0|\Theta_i)
  \1_{|Y_0|\wt\Theta_i\in A})_{i\ge 0})>1)-\P(c((g(|Y_0|\Theta_i)
  \1_{|Y_0|\wt\Theta_i\in A})_{i\ge 1})>1)
\big]
\Big|\to0\,.
\eeao
\ere
\subsection{ Beyond condition (C$_\vep$): Convergence of the tail empirical point process}\label{sec:tailemp}
In this section we consider an example, where the condition {\bf
  (C$_\vep$)}
in Theorem~\ref{th:antic} is not satisfied but Jakubowski's
telescoping sum approach is also applicable, yielding a
limit result. Related theory was developed in Balan and Louhichi 
\cite{balan:louhichi:2009} beyond the framework of \regvary\ \seq s,
in the context of triangular arrays of strictly stationary \seq
s and infinitely divisible limit laws; see also
Jakubowski and Kobus \cite{jakubowski:kobus:1989}.

Recycling notation, we 
define the random measures 
\beam\label{eq:process}
N_{nj}= k_n^{-1}\sum_{t=1}^j \vep_{a_m^{-1} X_t}\,,\quad 
j=1,\ldots,n,\quad N_{nn}=N_n\,,
\eeam
where $m=m_n\to\infty$ and $k_n=[n/m]\to\infty$.
The {\em tail empirical point process} $N_n$ plays an important role
in extreme value statistics; see Resnick and \sta\ 
\cite{resnick:starica:1998}, Resnick \cite{resnick:2007}, 
Drees and Rootz\'en \cite{drees:rootzen:2010} and the references therein.
\bth\label{thm:hill}
Consider  a strictly stationary $\bbr^d$-valued \seq\ $(X_t)$,  
satisfying the following conditions:
\begin{enumerate}
\item[\rm 1.]
${\bf (RV_\alpha)}$ for some $\alpha>0$, 
\item[\rm 2.]
the mixing  condition ${\mathcal A}(a_n)$ modified for the random \ms s
\eqref{eq:process},
\item[\rm 3.]
the anti-clustering
condition  \eqref{eq:antic} with $\Lambda_n=\{a_n\}$.
\end{enumerate}
Then the relation
$N_n\stp \mu_1$ holds in the space of random \ms s on $\ov\bbr_0^d$ 
equipped with the vague topology, where $n\,\P(a_n^{-1} X\in\cdot
)\stv \mu_1$ as $\nto$.
\ethe
The proof is given in Section~\ref{proofofthm:hill}.
\bre
Closely related results can be found in Resnick and \sta\
\cite{resnick:starica:1998} in the 1-dimensional case. 
Their mixing and anti-clustering
conditions are slightly different and they prove results for
triangular arrays under a vague tightness condition (which is
satisfied 
for  ${\bf (RV_\alpha)}$), similar to 
Balan and Louhichi \cite{balan:louhichi:2009}.  
It follows from the results in Resnick and \sta\ 
\cite{resnick:starica:1998}  that Theorem~\ref{thm:hill} implies the
consistency of the Hill estimator of $\alpha$ in the case  of positive
\rv s $(X_t)$, i.e., if $m_n\to\infty$ and $m_n/n\to 0$ then
\beao
\Big(\sum_{t=1}^{m_n-1} \log \Big(X_{(n-i+1)}/X_{(n-m_n+1)}\Big)\Big)^{-1}\stp
\alpha\,,\quad \nto\,,
\eeao
where $X_{(1)}\le \cdots \le X_{(n)}$ a.s. is the ordered sample of $X_1,\ldots,X_n$.
\ere

\section{Proofs}\label{sec:proofs}
\subsection{Proof of Theorem~\ref{th:antic}}\label{sec:proofofthm1}
For a given integer $n\ge 1$ and $\delta>0$, we define  
\beao
c_x(j_1,j_2)&=&\left\{\barr{ll}f(x^{-1}X_{j_1},\ldots,x^{-1}X_{j_2})&  
1\le j_1\le j_2\le  n\\ 
0 & j_1>j_2,\earr\right.\\
h_k(x^{-1}X_j,\ldots,x^{-1}X_n)&=&\left\{\barr{ll}
\1_{|X_j|\le x\delta} & k=0\\
\1_{M_{j+k,n}\le x\delta} &k\ge 1\,.\earr\right.\,,\quad j\le n\,. \eeao
By a telescoping sum argument, we decompose
\beao
c_x(1,n)=\sum_{j=1}^{n}[c_x(j,n)-c_x(j+1,n)]\,.
\eeao
We denote $\Delta_x(n)=c_x(1,n)-\sum_{j=1}^n[c_x(j,j+k)-c_x(j+1,j+k)]$ for some $k\ge 1$. By stationarity, we have
\beao
\E \Delta_x(n)= 
\E f(x^{-1}X_{1},\ldots,x^{-1}X_{n})- n\,
\big[\E f(x^{-1}X_{0},\ldots,x^{-1}X_{k})-\E f(x^{-1}X_{1},\ldots,x^{-1}X_{k})\big]\,.
\eeao
Moreover, as $f$ satisfies {\bf
  (C$_\vep$)} and  assuming without loss of
  generality that $|f|\le 1$, we obtain for any $\delta\le \vep$ that
\beao
|\Delta_x(n)|&=&\Big| \sum_{j=1}^{n}[c_x(j,n)-c_x(j+1,n)]-[c_x(j,j+k)-c_x(j+1,j+k)]\Big|\\
&\le&\Big|\sum_{j=1}^{n}c_x(j,n)(1-h_0(x^{-1}X_j,\ldots,x^{-1}X_n))(1-h_{k}(x^{-1}X_j,\ldots,x^{-1}X_n))\Big|\\
&\le&\sum_{j=1}^n\ \1_{|X_j|> x\delta} \1_{M_{j+k,n}> x\delta}. 
\eeao
Using stationarity,  we obtain
\beao
\E |\Delta_x(n)|&\le &  
\sum_{j=1}^n\P(|X_j|>x\delta\,, M_{j+k,n}>x\delta)\\
&\le &  n\,\P(|X_0|>x\delta\,,M_{k,n}>x\delta)\,.
\eeao
Using the uniform anti-clustering condition \eqref{eq:antic}, we
conclude that
\beao
\lim_{\kto}\lim_{\nto}\sup_{x\in\Lambda_n}\dfrac{\E
  |\Delta_x(n)|}{n\,\P(|X_0|>x)}=0\,,\quad \delta>0\,.
\eeao
It remains to prove the existence of the limit, uniformly for $x\in\Lambda_n$,
\beao
\lim_{k\to\infty}\lim_{n\to\infty}\frac{\E [f(x^{-1}X_{0},\ldots,x^{-1}X_{k})-
f(x^{-1}X_{1},\ldots,x^{-1}X_{k})]}{\P(|X_0|>x)},
\eeao
and to identify it.
By \regvar\ and  {\bf (C$_\vep$)}, we have
\beao
\lefteqn{\sup_{x\in \Lambda_n}
\Big|\dfrac{\E [f(x^{-1}X_{0},\ldots,x^{-1}X_{k})-
f(x^{-1}X_{1},\ldots,x^{-1}X_{k})]}{\P(|X_0|>x)}}\\
&&-\int_{0}^\infty\E [f(y\Theta_{0},\ldots,y\Theta_{k})
- f(0,y\Theta_1,\ldots,y\Theta_{k})  ]\,d(-y^{-\alpha})\Big|\\
&=&\lefteqn{\sup_{x\in \Lambda_n}
\Big|\dfrac{\E \Big[\big(f(x^{-1}X_{0},\ldots,x^{-1}X_{k})-
f(x^{-1}X_{1},\ldots,x^{-1}X_{k})\big)
\1_{|X_0|>x}\Big]}{\P(|X_0|>x)}}\\
&&-\int_{0}^\infty\E [f(y\Theta_{0},\ldots,y\Theta_{k})
- f(0,y\Theta_1,\ldots,y\Theta_{k})  ]\,d(-y^{-\alpha})\Big|\\
&\to&
0\,,\quad \nto\,.
\eeao
Finally, for every $k\ge 1$, in view of  {\bf (C$_\vep$)},
\beam\label{eq:dobbel}
\lefteqn{\int_{0}^\infty\E [f(y\Theta_{0},\ldots,y\Theta_{k})
- f(0,y\Theta_1,\ldots,y\Theta_{k})  ]\,d(-y^{-\alpha})}\nonumber\\&=&
\int_{\vep}^\infty\E [f(y\Theta_{0},\ldots,y\Theta_{k})
- f(0,y\Theta_1,\ldots,y\Theta_{k})  ]\,d(-y^{-\alpha})\,.
\eeam
The absolute value of the integrand is bounded by 2, hence integrable on
$[\vep,\infty)$. Moreover, under the anti-clustering condition
\eqref{eq:antic}  we have 
$\Theta_k\stp0$ as $k\to\infty$ as follows from the Remark \ref{rem:ac+} above.
Therefore and by {\bf (C$_\vep$)} the limits 
$\lim_{\kto}f(y\Theta_0,\ldots,y \Theta_{k})$ exist
and are finite for $y>0$. Dominated \con\ implies 
that one may let $k\to\infty$ in \eqref{eq:dobbel} and interchange the
limit and the integral.
\par
Combining the  partial limit results above, we conclude that the
theorem is proved. 
\subsection{Proof of Theorem~\ref{th:kdep}}\label{sec:proofofth2}
We start by 
collecting some properties of the \seq\ $(\Theta_t)$ which 
are specific for a $k$-dependent \regvary\ sequence. In particular, 
we will show that the conditional \pro y laws \eqref{eq:Theta} 
are well defined. 
\bpr\label{lem:1}
Assume the stationary $\bbr^d$-valued $k$-dependent  \seq\ $(X_t)$  
is such that $(X_0,\ldots,X_k)$ is \regvary\ with index $\alpha$ then $(X_t)$  satisfies  ${\bf (RV_\alpha)}$. The following properties also hold: 
\begin{itemize}
\item $\P(\Theta_t=0)=1$ for $|t|\ge k+1$,
\item For $1\le |t|\le k$, $\P(\Theta_t\neq0,\Theta_{t+j}\neq 0)=0$ for $|j|\ge k+1$,
\item  $\P(\Theta_{-t}=0,t=1,\ldots,k)>0$.
\end{itemize}
\epr
\begin{proof}[Proof of Proposition~\ref{lem:1}] Assume that $(X_0,\ldots,X_k)$ is \regvary\ and that $(X_t)$ is $k$-dependent. Then, for any $|t|>k$, we have by independence of $X_t$ and $X_0$ that $\P(\Theta_t=0)=1$. As the limit law is degenerate, the convergence in the definition of $\Theta_t$ holds also in probability; $\P(|X_t|/|X_0|\le\vep \mid |X_0|>x)\to 1$ for $\vep>0$.  By an application of a Slutsky argument, as $(X_0,\ldots,X_k)/|X_0|$ converges in distribution, conditionally on $|X_0|>x$, then it is also true that $(X_0,\ldots,X_t)/|X_0|$ given  $|X_0|>x$ converges to $(\Theta_0,\ldots,\Theta_t)= (\Theta_0,\ldots,\Theta_k,0,\ldots,0)$.\\

The second property of the spectral tail process follows from the fact that
$$
\P(|Y_t|>\vep\,,|Y_{t-j}|> \vep )=0,\quad \vep>0\,,\quad
|j|\ge k+1\,.
$$
Indeed, by independence of $X_t$ and $X_{t-j}$ for $|j|>k$, 
\beao
\P(  |X_t|\wedge |X_{t-j}|>\vep a_n\mid |X_0|>a_n)&\le& \frac{\P(  |X_t|\wedge |X_{t-j}|>\vep a_n)}{\P(|X_0|>a_n)}\\
&=& \frac{[\P(  |X_0|>\vep a_n)]^2}{\P(|X_0|>a_n)}\sim \vep^{-2\alpha}
\P(|X|>a_n)\to 0\,.
\eeao
The second property implies that
$$
\P(\max_{-k\le  t<0}|\Theta_t|>0,\Theta_{k}\neq 0)=0.
$$
Then
\beao
\P(\Theta_{k}\neq 0)=\P(\max_{-k\le t<0}|\Theta_t|=0,\Theta_{k}\neq 0)\le  \P(\max_{-k\le t<0}|\Theta_t|=0)\,,
\eeao
and the third property follows if $\P(\Theta_{k}\ne 0)>0$. 
Now assume that $\P(\Theta_{k}\neq 0)=0$. 
By the time change formula in Basrak and Segers \cite{basrak:segers:2009},
$\P(\Theta_{k}\neq
0)=\E|\Theta_{-k}|^\alpha$ and $\P(\Theta_{-k}=0)=1$. Thus 
$\max_{-k\le t<0}|\Theta_t|=\max_{-k< t<0}|\Theta_t|$ a.s. A recursive
argument yields the third property in the general case.
\end{proof}

\begin{proof}[Proof of Theorem~\ref{th:kdep}] We apply
Theorem~\ref{th:antic} for $\Lambda_n=[x_n,\infty)$. By virtue of $k$-dependence 
the anti-clustering condition is trivially satisfied.
In view of Proposition~\ref{lem:1} it remains to show that 
\beam\label{eq:step2}
\lefteqn{\la_k\int_{0}^\infty\E(f_{k+1}(y\wt\Theta_0,\ldots,y\wt\Theta_{k}))d(-y^{-\alpha})}\nonumber\\&=&
\int_{0}^\infty\E[f_{k+1}(y\Theta_{0},\ldots,y\Theta_{k})
- f_{ k+1}(0,y\Theta_1,\ldots,y\Theta_{k})  ]d(-y^{-\alpha})\,.\eeam
We observe that
\beam\label{eq:vip}
\lefteqn{f_{k+1}(y\Theta_0,\ldots,y\Theta_{k})\1_{\Theta_{-j}=0,j=1,\ldots,k}}\nonumber\\
&=&f_{k+1}(y\Theta_0,\ldots,y\Theta_{k})-\sum_{t=1}^{k} f_{k+1}(y\Theta_0,\ldots,y\Theta_{k})\1_{\Theta_{-k+t}\neq0,\Theta_{-k+t+1}=0,\ldots,\Theta_{-1}=0}.
\eeam
From $k$-dependence we conclude that 
$\Theta_t\neq 0$ implies $\Theta_{t+j}=0$, $j\ge k+1$. Therefore
\beao
\lefteqn{f_{k+1}(y\Theta_0,\ldots,y\Theta_{k})\1_{\Theta_{-k+t}\neq0,\Theta_{-k+t+1}=0,\ldots,\Theta_{-1}=0}}\\&=&f_{k+1}(y\Theta_0,\ldots,y\Theta_{t},0,\ldots,0)\1_{\Theta_{-k+t}\neq0,\Theta_{-k+t+1}=0,\ldots,\Theta_{-1}=0}.
\eeao
By the time change formula in Theorem 3.1 (iii) in Basrak and Segers \cite{basrak:segers:2009} and the property {\bf (C$_\vep$)},
\beao
\lefteqn{\E[f_{k+1}(y\Theta_0,\ldots,y\Theta_{t},0,\ldots,0)\1_{\Theta_{-k+t}\neq0,\Theta_{-k+t+1}=0,\ldots,\Theta_{-1}=0}]}\\
&=&\E[f_{k+1}(y\Theta_{k-t}|\Theta_{k-t}|^{-1},\ldots,y\Theta_{k}|\Theta_{k-t}|^{-1},0,\ldots,0)|\Theta_{k-t}|^{\alpha}\1_{\Theta_{0}\neq0,\Theta_{1}=0,\ldots,\Theta_{k-t-1}=0}]\\
&=&\E[f_{k+1}(y\Theta_{k-t}|\Theta_{k-t}|^{-1},\ldots,y\Theta_{k}|\Theta_{k-t}|^{-1})|\Theta_{k-t}|^{\alpha}\1_{\Theta_{1}=0,\ldots,\Theta_{k-t-1}=0}]\,.
\eeao
Now, by Fubini's Theorem when $\Theta_{k-t}\neq0$, first integrating \wrt\  $y$ and then changing variables:
\beao
&&\int_{0}^\infty
\E
f_{k+1}(y\Theta_{k-t}|\Theta_{k-t}|^{-1},\ldots,y\Theta_{k}|\Theta_{k-t}|^{-1})\1_{\Theta_{1}=0,\ldots,\Theta_{k-t-1}=0,\Theta_{k-t}\neq0,}
d(-(y
|\Theta_{k-t}|^{-1})^{-\alpha})\\
&=&\int_{0}^\infty \E f_{k+1}(y\Theta_{k-t},\ldots,y\Theta_{k})\1_{\Theta_{1}=0,\ldots,\Theta_{k-t-1}=0,\Theta_{k-t}\neq0}\;d(-y^{-\alpha})\,.
\eeao
Taking into account \eqref{eq:vip} and the previous identities, we obtain
\beao
&&\hspace{-2cm}\int_{0}^\infty\E(f_{k+1}(y\Theta_0,\ldots,y\Theta_{k})\1_{\Theta_{-j}=0,j=1,\ldots,k})d(-y^{-\alpha})
\\
&=&\int_{0}^\infty \E[f_{k+1}(y\Theta_{0},\ldots,y\Theta_k)]d(-y^{-\alpha})\\
&&-\sum_{i=1}^{k}\int_{0}^\infty\E[f_{k+1-i}(y\Theta_{i},\ldots,y\Theta_k)\1_{\Theta_j=0,j=1,\ldots,i-1,\Theta_{i}\neq0 }]d(-y^{-\alpha})\,.
\eeao
Since $(f_\ell)$ satisfies the consistency property {\bf (C$_\vep$)} the \rhs\ turns into
$$
\int_{0}^\infty
\E[f_{k+1}(y\Theta_{0},\ldots,y\Theta_k)]d(-y^{-\alpha})
-\int_{0}^\infty\E\Big[f_ {k+1}(0,y\Theta_{1},\ldots,y\Theta_k)\sum_{i=1}^{k}\1_{\Theta_j=0,j=1,\ldots,i-1,\Theta_{i}\neq0}\Big]d(-y^{-\alpha})$$
and, as $f(0,\ldots,0)=0$, the desired result follows. \end{proof}
\subsection{Proof of Theorem~\ref{thm:alphaclt}}\label{sec:proofofth}
We start with the case $\alpha\in (0,1)$. In this case,
  the negligibility condition 
\beao
\lim_{\vep\downarrow 0}\limsup_{\nto} \P(|a_n^{-1} \ov
S_n|>\delta)=0\,,\quad \delta>0\,,
\eeao 
is satisfied. Indeed,
  an application of Markov's inequality and Karamata's theorem yields
\beao
\P(|a_n^{-1} \ov
S_n|>\delta)\le \dfrac{\E |X|\1_{|X|\le\vep a_n }}{a_n\vep \P(|X|>\vep
a_n)} \vep n\,\P(|X|>\vep a_n)\to c\, \vep^{1-\alpha}\,,\quad \nto\,,
\eeao
and the \rhs\ vanishes as $\vep\downarrow 0$.
Therefore we may focus on the limit behavior of the \seq\ 
$(a_n^{-1}\underline S_n)$. Fix a small value $\vep\in (0,1)$. The
mixing condition \eqref{eq:mix} implies that
\beao
\log \E \ex^{is' \underline S_n/a_n}\sim k_n \log \E \ex^{is'
  \underline S_m/a_n}\sim - k_n \Big(1- \E \ex^{is'
  \underline S_m/a_n}\Big)\sim \dfrac{\E \ex^{is'
  \underline S_m/a_n}-1}{m\,\P(|X|>a_n)}\,.
\eeao 
We define the \fct s 
\beao
f_l(x_1, \ldots,x_l)= \exp\Big(iy s'\sum_{t=1}^l \underline x_t\Big)-1\,,\quad l\ge 0\,,
\eeao
where  $\underline x_t = x_t \1_{|x_t|>\vep a_n}$. The \seq\
$(f_l)$ satisfies {\bf (C$_\vep$)} and $m \P(|X|>a_n)\to 0$ as $\nto$.
An application of Theorem~\ref{thm:alphaclt} yields
\beao
\dfrac{\E \ex^{is'
  \underline S_m/a_n}-1}{m\,\P(|X|>a_n)} \to \int_0^\infty
\E\Big[\ex^{iys' \sum_{j=0}^\infty \underline \Theta_j}-
 \ex^{iys' \sum_{j=1}^\infty \underline \Theta_j}
\Big]\,d(-y^{-\alpha})\,,\quad s\in\bbr^d\,,\quad \nto\,.
\eeao
Here $\underline \Theta_i=\Theta_i \1_{y |\Theta|>\vep}$.
To complete the proof we have to justify that we can let
$\vep\downarrow 0$ in the limiting expression. We observe that the
integrand vanishes on the event $\{|y\Theta_0|\le \vep\}=\{y\le
\vep\}$. 
Therefore 
the \rhs\ turns into
\beao
 \int_\vep^\infty
\E\Big[\big(\ex^{iys' \Theta_0}-1\big) \ex^{iys' \sum_{j=1}^\infty
  \underline \Theta_j}\Big] \,d(-y^{-\alpha})\,.
\eeao
The integrand is uniformly bounded and therefore integrable at
infinity. In order to apply dominated \con\
as $\vep\downarrow 0$, we observe that for some constant $c>0$,
\beao
 \int_\vep^1
\E\Big|\big(\ex^{iys' \Theta_0}-1\big) \ex^{iys' \sum_{j=1}^\infty
  \underline \Theta_j}\Big| \,d(-y^{-\alpha})\le  \int_\vep^1
\E |ys'\Theta_0|\,  d(-y^{-\alpha})\le c\,\int_0^1 y^{-\alpha}\, dy\,.
\eeao
Now an application of the dominated \con\ theorem as $\vep\downarrow
0$ yields the desired log-\chf\ of an $\alpha$-stable law.
\par
The proof in the case $\alpha\in (1,2)$ is similar. In view of the
negligibility condition \eqref{eq:neg}  we may focus on the limit behavior
of $(a_n^{-1}(\underline S_n- \E\underline S_n))$. Since $\E
\underline S_n/a_n$ converges,  
the mixing condition remains valid for the corresponding centered sums.
Then
\beao
\log \E \ex^{is'(\underline S_n-\E \underline S_n)/a_n}\sim 
\dfrac{\E \ex^{is (\underline S_m-\E \underline S_m)/a_n}-1}{m\,\P(|X>a_n)}\,.
\eeao
We also have 
\beao
\lefteqn{\Big|\Big(\E \ex^{is' (\underline S_m-\E \underline
    S_m)/a_n}-1\Big)- \Big(\E \ex^{is' \underline S_m/a_n}-1 - is' \E
  \underline S_m/a_n\Big)\Big|}\\&=&
\Big|\Big(\E \ex^{is' \underline S_m/a_n}-1\Big) \Big(\ex^{-is'\E \underline
    S_m/a_n}-1\Big) + \Big(\ex^{-is'\E \underline
    S_m/a_n}-1+is' \E
  \underline S_m/a_n\Big)\Big|\\
&\le & c\big(\E|\underline S_m |/a_n\big)^2\\
&\le & c \big(m \E|X/a_n|\1_{|X|>\vep
  a_n}\big)^2=O(m\,\P(|X|>a_n))=o(1)\,,\quad \nto\,.
\eeao
and  
\beam\label{eq:expect}
\dfrac{a_n^{-1}s'\E \underline S_m}{m\,\P(|X|>a_n)}&=&\dfrac{s'\E \underline
  X_0}{a_n \P(|X|>a_n)}\nonumber\\
&\sim & \E (s'X_0/(\vep a_n)\mid |X_0|>\vep a_n) \vep^{1-\alpha}\nonumber\\
&\to & s'\E Y_0\,\vep^{1-\alpha}\nonumber\\
&=& \E |Y_0|\,s'\E \Theta_0 \,\vep^{1-\alpha}\nonumber\\
&=& s'\E \Theta_0\,\dfrac{\alpha}{\alpha-1}\vep^{1-\alpha}\,.
\eeam
An application of Theorem~\ref{th:antic} yields as $\nto$,
\beao\lefteqn{
\dfrac{\E\Big[ \ex^{is'
  \underline S_m/a_n}-1- is'\E \underline S_m/a_n\Big]}{m\,\P(|X|>a_n)}}\\& \to&
\int_0^\infty \,
\E\Big[\Big(\ex^{iys' \sum_{j=0}^\infty
  \underline \Theta_j} - 1- iys' \sum_{j=0}^\infty
  \underline \Theta_j\Big)- 
\Big(\ex^{iys' \sum_{j=1}^\infty
  \underline \Theta_j } - 1- iys' \sum_{j=1}^\infty
  \underline \Theta_j\Big)\Big]\,
d(-y^{-\alpha})\\
&=&\int_\vep^\infty \,
\E\Big[\Big(\ex^{iys' \Theta_0}-1\Big)\Big(\ex^{iys' \sum_{j=1}^\infty
  \underline \Theta_j}-1\Big) + \Big(\ex^{iys' \Theta_0}-1-iys' 
  \Theta_0\Big)\,\Big]
d(-y^{-\alpha})\,.
\eeao
The integrand is bounded by $cy$ for large values of $y$ and therefore
integrable at infinity. We also have 
\beao
&&\int_\vep^1 \,
\E\Big|\Big(\ex^{iys' \Theta_0}-1\Big)\Big(\ex^{iys' \sum_{j=1}^\infty
  \underline \Theta_j}-1\Big) + \Big(\ex^{iys' \Theta_0}-1-iys' 
  \Theta_0\Big)\,\Big|
d(-y^{-\alpha})\\
&\le &c\,\Big(\sum_{j=1}^\infty
  \E|\Theta_j|+1\Big)\,\int_0^1 y^2 \,d(-y^{-\alpha})<\infty
\eeao
Therefore an application of the
dominated \con\ theorem yields the desired log-\chf\ of a stable law
in the case $\alpha\in (1,2)$.

\subsection{Proof of Theorem~\ref{thm:suprema}}\label{prof}
We start with the case $\alpha>1$.
We have for any $\delta>0$, 
\beam
 \lefteqn{\P\Big(\sup_{t\le n} \ul S_t>(1+\delta) x\Big)- \P\Big(\sup_{t\le n} |\ov S_t|>\delta x\Big)}\nonumber\\ &\le& \P\Big(\sup_{t\le n}  S_t>  x\Big)\le    \P\Big(\sup_{t\le n} \ul S_t>(1-\delta) x\Big)+\P\Big(\sup_{t\le n} |\ov S_t|>\delta x\Big).\label{eq:99}
\eeam
In view of condition  \eqref{eq:nega},  the limiting behavior will be determined by the ratios
$$
\frac{ \P(\sup_{t\le n } \ul S_t>(1\pm
    \delta)x)}{n\P(|X_0|>x)}=\frac{ \E(f_n(x^{-1}X_1,\ldots,x^{-1}X_n))}{n\P(|X_0|>x)}
$$
with the \fct s $f_l(x_1,\ldots,x_l)= \1_{\sup_{t\le l}
  (\underline x_1+\cdots +\underline x_t)>(1\pm \delta)}$ for fixed
small 
$\delta,\vep>0$. The functions $f_l$ satisfy the condition (C$_\vep$) and 
an application of Theorem~\ref{th:antic} yields for
$\vep<1$ as $\nto$,
\beao
\lefteqn{\sup_{x\in \Lambda_n} \Big|\dfrac{\P\big(\sup_{i\le n}
  \underline S_i>(1\pm \delta)x\big)}{n\,\P(|X|>x)}}\\&& - (1\pm\delta)\int_0^\infty\Big[
\P\Big(\Theta_0+\sup_{t\ge 1} \sum_{i=1}^t
\underline\Theta_i>y^{-1}\Big)-
\P\Big(\sup_{t\ge 1} \sum_{i=1}^t
\underline\Theta_i>y^{-1}\Big)\Big]\, d(-y^{-\alpha})
\Big|\to 0\,.
\eeao
Finally, we can take limits as $\vep,\delta \downarrow 0$, using a
domination argument. The domination argument is justified because
we have 
\beao
\lefteqn{\int_0^\infty\Big[
\P\Big(\Theta_0+\sup_{t\ge 1} \sum_{i=1}^t
\underline\Theta_i>y^{-1}\Big)-
\P\Big(\sup_{t\ge 1} \sum_{i=1}^t
\underline\Theta_i>y^{-1}\Big)\Big]\, d(-y^{-\alpha})}\\
&=&
\alpha\int_0^\infty\Big[
\P\Big(\Theta_0+\sup_{t\ge 1} \sum_{i=1}^t
\underline\Theta_i>z\Big)-
\P\Big(\sup_{t\ge 1} \sum_{i=1}^t
\underline\Theta_i>z\Big)\Big]\, z^{\alpha-1}\,dz \\
&=&\alpha\int_0^\infty
\E\Big[\Big(\1_{\Theta_0=1,1+\sup_{t\ge 1} \sum_{i=1}^t
\underline\Theta_i>z}-
\1_{\Theta_0=1,\sup_{t\ge 1} \sum_{i=1}^t
\underline\Theta_i>z}\Big)\\
&&+\Big(\1_{\Theta_0=-1,-1+\sup_{t\ge 1} \sum_{i=1}^t
\underline\Theta_i>z}-
\1_{\Theta_0=-1,\sup_{t\ge 1} \sum_{i=1}^t
\underline\Theta_i>z}\Big)\Big]\, z^{\alpha-1}\,dz \\
&=& \E\Big[\1_{\Theta_0=1} \Big[\Big(1+\sup_{t\ge 1} \sum_{i=1}^t
\underline\Theta_i\Big)_+^\alpha-\Big(\sup_{t\ge 1} \sum_{i=1}^t
\underline\Theta_i\Big)_+^\alpha\Big]\\
&&- \1_{\Theta_0=-1} \Big[\Big(\sup_{t\ge 1} \sum_{i=1}^t
\underline\Theta_i\Big)_+^\alpha-\Big(-1+\sup_{t\ge 1} \sum_{i=1}^t
\underline\Theta_i\Big)_+^\alpha\Big]\Big]\\
&=& \E\Big[\Big(\Theta_0+\sup_{t\ge 1} \sum_{i=1}^t
\underline\Theta_i\Big)_+^\alpha-\Big(\sup_{t\ge 1} \sum_{i=1}^t
\underline\Theta_i\Big)_+^\alpha\Big]\,.
\eeao
Also observe that 
\beam\label{eq:al0}
\lefteqn{\E\Big|\Big(\Theta_0+\sup_{t\ge 1} \sum_{i=1}^t
\underline\Theta_i\Big)_+^\alpha-\Big(\sup_{t\ge 1} \sum_{i=1}^t
\underline\Theta_i\Big)_+^\alpha\Big|}\\&\le& c\Big[1+ \E \Big(\sum_{i=1}^\infty
|\un \Theta_i|\Big)^{\alpha-1}\Big]\le c\,\Big[1+ \E \Big(\sum_{i=1}^\infty
|\Theta_i|\Big)^{\alpha-1}\Big]\,. \nonumber
\eeam
Under the assumptions of the theorem, the \rhs\ is finite. An
application of Lebesgue dominate \con\ yields that
\beao
\E\Big[\Big(\Theta_0+\sup_{t\ge 1} \sum_{i=1}^t
\underline\Theta_i\Big)_+^\alpha-\Big(\sup_{t\ge 1} \sum_{i=1}^t
\underline\Theta_i\Big)_+^\alpha\Big]\to \E\Big[\Big(\Theta_0+\sup_{t\ge 1} \sum_{i=1}^t
\Theta_i\Big)_+^\alpha-\Big(\sup_{t\ge 1} \sum_{i=1}^t
\Theta_i\Big)_+^\alpha\Big]\,,\quad \vep \to 0\,.
\eeao
This concludes the proof in the case $\alpha>1$.
\par
In the case $\alpha\in (0,1]$
one can
follow the lines of the proof but instead of \eqref{eq:al0} one can use
concavity to obtain the bound
\beao
\E\Big|\Big(\Theta_0+\sup_{t\ge 1} \sum_{i=1}^t
\underline\Theta_i\Big)_+^\alpha-\Big(\sup_{t\ge 1} \sum_{i=1}^t
\underline\Theta_i\Big)_+^\alpha\Big| \le c\,.
\eeao
An application of Lebesgue dominated \con\ finishes  the proof.

\subsection{Proof of Corollary~\ref{th:sup sum}}\label{proofofcorol}
We apply Theorem~\ref{thm:suprema}. For a $k$-dependent \seq , 
the anti-clustering condition  \eqref{eq:antic} is trivially satisfied
and we also have $\E |\Theta_i|^\alpha<\infty$ for all $i$. 
The vanishing-small-values condition is verified in Lemma~\ref{lem:8} below.
Now one can follow the lines of the proof of Theorem~\ref{thm:suprema} but
instead of Theorem~\ref{th:antic} apply Theorem~\ref{th:kdep} to
obtain
\beao
\lefteqn{\sup_{x\ge x_n}\Big|\frac{ \P(\sup_{t\le n } \ul S_t>(1\pm
    \delta)x)}{n\P(|X_0|>x)}}\\&&- (1\pm \delta)^{-\alpha}\lambda_k \int_0^\infty \P\Big( \sup_{t\le k}
y \sum_{i=0}^t\Theta_i \1_{y |\Theta_i|>\vep } >1\Big) d(-y^{-\alpha})\Big|\to 0.
\eeao
Finally, we can take limits as $\vep,\delta \downarrow 0$, using a
domination argument and the fact that $\E|\Theta_i|^\alpha<\infty$ for
all $i$.
\ble\label{lem:8} Assume the  conditions of Corollary~\ref{th:sup sum}.
Then
\beam\label{eq:svv}
\lim_{\vep\downarrow 0}\limsup_{\nto}\sup_{x\ge x_n} \frac{\P\Big(\sup_{t\le n} |\ov S_t|> x\Big)}{n\P(|X|>x)}=0\,.
\eeam
\ele
\begin{proof}[Proof of Lemma~\ref{lem:8}] 
We observe that
$$
\sup_{t\le n} |\ov S_t|\le \sum_{j=1}^{k+1}\sup_{s(k+1)\le n}|\sum_{t=0}^s\ov X_{(k+1)t+j}|\,.
$$
Therefore for $\delta>0$,
\beao
\frac{\P\Big(\sup_{t\le n} |\ov S_t|>\delta x\Big)}{n\P(|X_0|>x)}
&\le & \frac{(k+1)  \P\Big(\sup_{s(k+1)\le n}|\sum_{t=0}^s\ov X_{(k+1)t}|>\delta
  x/(k+1)\Big)}{n\P(|X_0|>x)}=  I(x)\,.
\eeao
If $\alpha\in (0,1)$, we have by Markov's inequality,
\beao
\sup_{x\ge x_n}I(x)&\le&\sup_{x\ge x_n} \frac{(k+1) \P\Big(\sum_{t=1}^{[n/(k+1)]}|\ov X_{kt}|>\delta x/(k+1)\Big)}{n\P(|X_0|>x)}\\&\le&
 c   \sup_{x\ge x_n} \frac{ \E |\ov X_0|}{x\P(|X_0|>x)}\\
&\le & c  \sup_{x\ge x_n} \frac{ \E(|X_0|\1_{|X_0|\le x\vep})}{
   \E(|X_0|\1_{|X_0|\le x})}\sup_{x\ge x_n} \frac{ \E(|X_0|\1_{|X_0|\le
     x})}{x\P(|X_0|>x)}\,.
\eeao
The second supremum is bounded in view of Karamata's
theorem and the first supremum is bounded by  $c\vep^{1-\alpha}$
uniformly for small $\vep$ is view of the uniform \con\ theorem for
\regvary\ \fct s. Hence the \rhs\ converges to zero by first letting
$\nto$ and then $\vep\downarrow 0$.
\par
If $\alpha\ge 1$, we have 
\beao
I(x)\le  \frac{(k+1) \P\Big(\sup_{s(k+1)\le n}|\sum_{t=0}^s(\ov X_{(k+1)t}- \E \ov X)|>\delta
  (x- \delta^{-1} n x^{-1}|\E\ov X|)/(k+1)\Big)}{n\P(|X_0|>x)}\,.
\eeao
If 
$\E|X|<\infty$ and $\E X=0$ we have by Karamata's theorem and uniform \con\
\beao
n x^{-1}|\E\ov X|\le n x^{-1} \E|X| \1_{|X|>\vep x}\le c (\vep) 
n\P(|X|>x)=o(1)\,,\quad \nto\,.
\eeao
The last identity follows from the fact that $n\P(|X|>x_n)=o(1)$ as
$\nto$. In the case $\alpha=1$ and $\E |X|=\infty$ we require the
corresponding condition \eqref{eq:auxo}.
Therefore we have 
\beao
I(x)\le  \frac{(k+1)  \P\Big(\sup_{s(k+1)\le n }|\sum_{t=0}^s(\ov X_{(k+1)t}- \E
  \ov X)|>0.5 \delta
  x /(k+1)\Big)}{n\P(|X_0|>x)}\,.
\eeao
uniformly for $x\ge x_n$ and large $n$. To ease notation, we will
assume without loss of generality that $E\ov X=0$.
Now an application of the 
Fuk-Nagaev inequality (see Petrov \cite{petrov:1995}, p. 78) for
$p>\alpha\vee 2$ yields,
\beao
I(x)\le \frac{c \Big( n x^{-p} \E|\ov X |^p+ \ex^{-c x^2(n\E\ov X^2)^{-1}}\Big)}{n\P(|X|>x)}\le c
\frac{\E|\ov X |^p}{x^p\P(|X|>x)}+ c  \frac{\ex^{-c x^2(n\E\ov X^2)^{-1}}}{n\P(|X|>x)}\,.
\eeao
The first summand on the \rhs\ converges to $\vep^{p-\alpha}$
uniformly for small $\vep$ and uniformly for $x\ge x_n$. 
The second term  is uniformly negligible for $x\ge x_n$  
if $\alpha<2$ in view of the relation  $x^2 (n\E\ov
X^2)^{-1}\sim c(n\P(|X|>x))^{-1}$. If $\alpha=2$ and $\E X^2=\infty$,
Karamata's theorem yields that $\E (\ov X/x)^2/\P(|X|>x)$ is a
\slvary\ \fct\ converging to infinity as $\xto$, and then the growth
condition $x_n/n^{0.5+\delta}\to\infty$ ensures that the second term
is negligible. 
If $\E|X|^2<\infty$, 
writing $x=\sqrt{n\log n}y$ for $y>0$ sufficiently large, we obtain
$$
 \frac{\ex^{-c x^2(n\E\ov X^2)^{-1}}}{n\P(|X|>x)}\le n^{\alpha +\vep
   -cy^2-1}y^{\alpha+\vep}\le n^{-\eta}$$ for some $\eta>0$ and
 thus the \rhs\  converges to $0$ uniformly for $x\ge x_n$.
\end{proof}
\subsection{Proof of Theorem~\ref{thm:ruina}}\label{proofofruin}
We start by proving that
\beam\label{eq:r35}
\lim_{C\to\infty}\limsup_{\xto} 
\dfrac{\P(\sup_{t>Cx}(S_t-\rho t)>x)}{x\,\P(|X|>x)}&=&0\,.
\eeam
Assume that $[Cx]\in D_k=(2^{k},2^{k+1}]$ for some $k\ge 1$. Then
\beao
\P(\sup_{t>Cx}(S_t-\rho t)>x)&\le& \sum_{l=k}^\infty \P\big(\sup_{t\in
  D_l} S_t>x+\rho\,2^l\big)\\
&\le & \sum_{l=k}^\infty \P\big(\sup_{t\le
  2^{l+1}} S_t>x+\rho\,2^l\big)\,.
\eeao
We observe that $x+\rho\,2^l \in D_l$ and therefore we are in the
range where we may apply the \ld\ results for suprema in
Theorem~\ref{thm:suprema}. Then the \rhs\ above can be bounded
uniformly by
\beao
c\,\sum_{l=k}^\infty 2^l\,\P(|X|>x+\rho \,2^l)&\le& c \int_{Cx}^\infty
\P(|X|>\rho y +x)\,dy\\  
&\le & c\, (Cx)  \P(|X|>Cx)\sim \, c\,C^{1-\alpha}\,x\,\P(|X|>x)\,.
\eeao
In the last step we used Karamata's theorem.
Relation \eqref{eq:r35} follows by observing that $\alpha>1$. 
\par
In view of \eqref{eq:r35} it suffices to bound the \pro ies 
$\P(\sup_{t\le Cx}(S_t-\rho t)>x)$ as $\xto$ for any large value of
$C$. 
Moreover, in view of condition \eqref{eq:nega}
we may replace the random
walk $(S_t)$ in the ruin \pro y by the truncated version $(\underline
S_t)$, letting $\vep\downarrow 0$ in the final bound. We then adapt 
the proof of Theorem \ref{th:antic} to the case where the 
functional depends on the last index larger than $\vep>0$. 
More specifically, denoting
$$
c_{x,j_2}(j_1,j_2)=1_{\sup_{j_1\le t\le j_2}(\sum_{i=j_1}^{t}\underline X_i-\rho t)>x},\quad 1\le j_1\le j_2\le [Cx]
$$
we follow the reasoning of Section \ref{sec:proofofthm1}. Observing that
\beao
\E\big([c_{x,[Cx]}(j,[Cx])-c_{x,[Cx]}(j+1,[Cx])]-[c_{x,j+k}(j,j+k)-c_{x,j+k}(j+1,j+k)]\big)\\
\le \P(|X_j|> x\vep, M_{j+k,[Cx]}> x\vep),
\eeao
we use the stationarity and the anti-clustering condition
\eqref{eq:antic}  to estimate the bound by 
\beam\label{eq:k9}
\P( M_{k,[Cx]}>\vep x, |X_0|>\vep x)\le C_k\,\P(|X|>x) 
\eeam
for large $x$ and a real \seq\ $(C_k)$ (depending on $C$) \st\ $C_k\to
0$ as $\kto$.
Therefore it suffices to study the limiting behavior of the difference
of the two terms emerging from the telescoping argument when
$x\to\infty$ and then $k\to\infty$; we indicate how we treat the first one:
$$
\sum_{j=1}^{[Cx]}\E(c_{x,j+k}(j,j+k))=\sum_{j=1}^{[Cx]}\P\Big(\sup_{j\le t\le j+k}\Big(\sum_{i=j}^{t}\underline X_i-\rho t\Big)>x\Big)=\sum_{j=1}^{[Cx]}\P\Big(\sup_{0\le t\le k}\Big(\sum_{i=0}^{t}\underline X_i-\rho (t+j)\Big)>x\Big).
$$
We have
\beao
\P(\sup_{0\le t\le k}(\ul S_{t }+\ul X_0 )>x+\rho k)\le\P(\sup_{0\le t\le k}(\ul S_{t }+\ul X_0-\rho (t+j))>x)\le \P(\sup_{0\le t\le k}(\ul S_{t }+\ul X_0)>x+\rho j).
\eeao 
We obtain the following sandwich bound
\beao
 x\int_{1/x}^{C}\P(\sup_{0\le t\le k}(\ul S_{t }+\ul
X_0)>(1+u\rho)x+\rho k)du\le\sum_{j=1}^{[Cx]} \P(\sup_{0\le t\le k}(\ul S_{t }+\ul X_0-\rho (t+j))>x)\\
 \le x\int_0^{C+x^{-1}}\P(\sup_{0\le t\le k}(\ul S_{t }+\ul
X_0)>(1+u\rho)x)du\,.
\eeao
By first letting $\xto$, using classical arguments from regular
variation 
theory,  in particular the uniform \con\ theorem, we obtain
\beao
\frac{\sum_{j=1}^{[Cx]} \P(\sup_{0\le t\le k}(\ul S_{t }+\ul X_0-\rho (t+j))>x)}{x\P(|X|>x)}du\sim \int_{0}^{C}\frac{\E(\sup_{t\le k}\sum_{i=0}^t\underline{\Theta}_i)_+^\alpha}{(1+u\rho)^\alpha}du.
\eeao
We then determine the limit of the difference of the two terms
emerging 
from the telescoping argument. For any $k\ge1$ we have
$$
\frac{\P(\sup_{t\le Cx}(S_t-\rho t)>x)}{x\P(|X|>x)}=\int_{0}^{C}\frac{\E[(\sup_{t\le k}\sum_{i=0}^t\underline{\Theta}_i)_+^\alpha-(\sup_{1\le t\le k}\sum_{i=1}^t\underline{\Theta}_i)_+^\alpha]}{(1+u\rho)^\alpha}du+O(CC_k).
$$
The last term vanishes when $k\to\infty$ for any $C>0$. Moreover, one can use uniform convergence and let $\vep\to 0$.
Finally, letting $C\to \infty$ and observing that $\int_0^\infty
(1+\rho u)^{-\alpha}\,du= \rho^{-1} (\alpha-1)^{-1}$, we obtain the
desired ruin bound.

\subsection{Proof of Corollary~\ref{thm:ruin}}\label{proofofruina}
The proof follows the lines of the proof of  Theorem~\ref{thm:ruina}.
We mention that condition~\ref{eq:antic} is trivially satisfied and the
vanishing-small-values condition holds in view of Lemma~\ref{lem:8}. 
\subsection{Proof of Theorem~\ref{thm:hill}}\label{proofofthm:hill}
As in the proof of Theorem~\ref{cor:ppgeneral}, we
  have for any $g\in\bbc_K^+$ with $g(x)=0$ for $|x|\le \delta$ for
  some $\delta>0$, in view of the 
mixing condition 
${\mathcal A}(a_n)$, 
\beao
-\log \E \ex^{-\int g dN_n}\sim k_n \big(1-\E \ex^{-\int g
  dN_{nm}}\big)\,,\quad \nto\,.
\eeao 
Now we can proceed as in the beginning of the proof of Theorem~\ref{th:antic}.
For $k\ge 2$,
\beao
&&\Big|k_n \E \big(1-\ex^{-\int g dN_{nm}}\big) - 
n \big[\E\big(1-\ex^{-\int gdN_{n,k-1}}\big)- \E\big(1-\ex^{-\int
  gdN_{n,k}}\big)\big]\Big|\\&\le&
 \P(M_{k,m}>\delta a_m\mid |X_0|>a_m\delta)\,.
\eeao
The \rhs\ is negligible by virtue of \eqref{eq:antic}.
Applying a Taylor
expansion, we have for some random $\xi\in (0,1)$,
\beao
n \,\E\big(1-\ex^{-\int g dN_{nk}}\big)
&=&n\Big[ \E \int g dN_{nk}-0.5\E\Big[\ex^{-\xi \int g dN_{nk}}\Big( \int
g\,dN_{nk}\Big)^2\Big]\Big]=I_1-I_2\,.
\eeao
Let 
$\mu_{0,l}$, $l\ge 1$, be the limit \ms s of \regvar\ for $a_n^{-1}(X_0,X_l)$, then we have
\beao
I_1&=& m_n \,k \,\E g(a_m^{-1}X)\to k \,\int g\, d\mu_1\,,\\
I_2&\le & n k_n^{-2} \E\Big( \sum_{t=1}^k g(a_m^{-1} X_t)\Big)^2\\
&=& k_n^{-1} \Big[k \,m_n\,\E g^2(a_m^{-1} X)+ 2\, \sum_{h=1}^{k-1}
(k-h)\, m_n\, \E[ g(a_m^{-1} X_0) g(a_m^{-1} X_h)] \Big]\\
&\sim &  k_n^{-1} \Big[ k \int g^2\, d\mu_1 +  2 \sum_{h=1}^{k-1}
(k-h) \int g(x)\, g(y)\, \mu_{0,h} (dx,dy)\Big]
\to 0\,,\quad \nto\,.
\eeao
Finally, we have
 \beao
\lim_{\nto} k_n \,\E(1- \ex^{-\int g\, dN_{nm}})& =&\lim_{k\to \infty} 
\lim_{\nto}
n \,\Big[\E\big(1-\ex^{-\int gdN_{n,k-1}}\big)- \E\big(1-\ex^{-\int
  fdN_{n,k}}\big)\Big]\\& =&\int g \,d\mu_1\,.
\eeao
This proves the \con\ of the Laplace \fct als
\beao
\E \ex^{-\int g\, dN_n }\to \ex^{- \int g\, d\mu_1 }\,,\quad g\in \bbc_K^+\,,
\eeao 
hence $N_n\stp \mu_1$.
\bre
The proof shows that the condition ${\bf (RV_\alpha)}$ is not really
needed in this case. It suffices that $X_0$ is \regvary\ and
$n\,\P(a_n^{-1} (X_0,X_h)\in \cdot)\stv \mu_{0,h}$ for every $h\ge 1$,
where  $\mu_{0,h}$ is a Radon \ms\ on $\ov\bbr_0^d$ which, possibly, is
the null \ms .

\ere

{\noindent \bf Acknowledgement.} We would like to thank the referee for careful reading of the paper and for various suggestions which led to an improvement of the presentation.

\bibliographystyle{imsart-nameyear}

\end{document}